\input amssym.def
\input amssym.tex
% 09 / 01 / 07

\def\item#1{\vskip1.3pt\hang\textindent {\rm #1}}% THIS REPLACES KNUTH'S DEF'N

                                  % THIS REPLACES KNUTH'S DEF'N

\tolerance=300
\pretolerance=200
\hfuzz=1pt
\vfuzz=1pt

% Print out with -x300 -y50

% bei unseren Artikeln uebliches Format:
\magnification=\magstep0

%\magnification=1100 
% gibt groesseres Schriftbild

% Offsetwerte fuer Ausdrucke in Erlangen, hoffset um 0.6 in groesser machen
\hoffset=0.6in
\voffset=0.8in

%\baselineskip                 
%\hsize=5.8 true in 
\hsize=5.8 true in 
\vsize=8.5 true in
%\baselineskip
\parindent=25pt
\mathsurround=1pt
\parskip=1pt plus .25pt minus .25pt
\normallineskiplimit=.99pt

\countdef\revised=100
\mathchardef\emptyset="001F % THIS REPLACES KNUTH'S DEFINITION
\chardef\ss="19
\def\3{\ss}
\def\anf{$\lower1.2ex\hbox{"}$}
\def\frac#1#2{{#1 \over #2}}
\def\>{>\!\!>}
\def\<{<\!\!<}

\def\ssarr{\hbox to 30pt{\rightarrowfill}}
\def\sarr{\hbox to 40pt{\rightarrowfill}}
\def\arr{\hbox to 60pt{\rightarrowfill}}

\def\larr{\hbox to 60pt{\leftarrowfill}}
\def\Arr{\hbox to 80pt{\rightarrowfill}}

{}

\def\Co{\mathop{\rm Co}\nolimits}

\def\Ad{\mathop{\rm Ad}\nolimits}

\def\Aut{\mathop{\rm Aut}\nolimits}

\def\det{\mathop{\rm det}\nolimits}

\def\Diff{\mathop{\rm Diff}\nolimits}

\def\End{\mathop{\rm End}\nolimits}

\def\Gl{\mathop{\rm Gl}\nolimits}

\def\Hom{\mathop{\rm Hom}\nolimits}%
\def\id{\mathop{\rm id}\nolimits} % USED FOR IDENTITY FUNCTION

% USED FOR IMAGINARY PART OF COMPLEX NUMBERS

\def\OO{{\rm O}}

% USED FOR REAL PART OF COMPLEX NUMBERS

%\def\reach{\mathop{\rm reach}\nolimits}
%\def\SAut{\mathop{\rm SAut}\nolimits}
%\def\SEnd{\mathop{\rm SEnd}\nolimits}

%\def\HSp{\mathop{\rm HSp}\nolimits}
\def\Sp{\mathop{\rm Sp}\nolimits}

\def\SU{\mathop{\rm SU}\nolimits}

\def\tr{\mathop{\rm tr}\nolimits}% USED FOR TRACE OF MATRIX
% USED FOR TRACE OF MATRIX
%\def\trdeg{\mathop{\rm trdeg}\nolimits} 

\def\OO{{\rm O}}

\def\0{{\bf 0}}
\def\1{{\bf 1}}

\def\g{{\frak g}}

\def\h{{\frak h}}

\def\m{{\frak m}}

\def\q{{\frak q}}

\def\C{{\Bbb C}}

\def\K{{\Bbb K}}

\def\N{{\Bbb N}} 
 
\def\P{{\Bbb P}} 
 
\def\R{{\Bbb R}} 
 
\def\Z{{\Bbb Z}} 

\def\:{\colon}  %8.5.92
\def\.{{\cdot}}
\def\|{\Vert}
\def\bsk{\bigskip}

\def\giantskip{\vskip2\bigskipamount}
\def\gsk{\giantskip}

\def\msk{\medskip}

\def\ssk{\smallskip}

\def\bbr{\bigbreak}
\def\giantbreak{\par \ifdim\lastskip<2\bigskipamount \removelastskip
         \penalty-400 \giantskip\fi}

\def\nin{\noindent}
\def\cen{\centerline}
\def\pagebreak{\vskip 0pt plus 0.0001fil\break}
\def\linebreak{\break}

\def\hat{\widehat}

\def\eps{\varepsilon}
\def\epsilon{\varepsilon}

\def\nin{\noindent}

\def\pder#1,#2,#3 { {\partial #1 \over \partial #2}(#3)}
\def\pde#1,#2 { {\partial #1 \over \partial #2}}
\def\phi{\varphi}

% Besser \ltimes und \rtimes aus dem AMS-Symbols verwenden. 
%\def\sdir#1{\hbox{$\mathrel\times{\hskip -4.6pt 
%            {\vrule height 4.7 pt depth .5 pt}}\hskip 2pt_{#1}$}}

\def\tilde{\widetilde}

\font\eightrm=cmr8

% SANS SERIF 10 POINT
 %SANS SERIF 10 POINT ITALIC

%\font\smc8=cmcsc8 
 %SLANTED TYPEWRITER 10 POINT
 %BOLD FACE MATH SYMBOLS 10 POINT
 %DUNHILL STYLE 10 POINT
 %SAN SERIF BOLD EXTENDED 10 POINT
 %USED FOR TITLES
 %USED FOR TITLES
\font\bfone=cmbx10 scaled\magstep1 %BOLDFACE AT MAGSTEP 1
\font\bftwo=cmbx10 scaled\magstep2 %BOLDFACE AT MAGSTEP 2
 %BOLDFACE AT MAGSTEP 3

\def\qed{{\unskip\nobreak\hfil\penalty50\hskip .001pt \hbox{}\nobreak\hfil
          \vrule height 1.2ex width 1.1ex depth -.1ex
           \parfillskip=0pt\finalhyphendemerits=0\medbreak}\rm}
%This is the end-of-proof sign. 
%Not to be used in display mode. 
%If you want to conclude a proof 
%at the end of a line is display mode use 

\def\qeddis{\eqno{\vrule height 1.2ex width 1.1ex depth -.1ex} $$
                   \medbreak\rm}
%BUT OMIT $$---the macro will write that

\def\Lemma #1. {\bigbreak\vskip-\parskip\noindent{\bf Lemma #1.}\quad\it}

\def\Sublemma #1. {\bigbreak\vskip-\parskip\noindent{\bf Sublemma #1.}\quad\it}

\def\Proposition #1. {\bigbreak\vskip-\parskip\noindent{\bf Proposition #1.}
\quad\it}

\def\Corollary #1. {\bigbreak\vskip-\parskip\nin{\bf Corollary #1.}
\quad\it}

\def\Theorem #1. {\bigbreak\vskip-\parskip\noindent{\bf Theorem #1.}
\quad\it}

\def\Definition #1. {\rm\bigbreak\vskip-\parskip\noindent{\bf Definition #1.}
\quad}

\def\Remark #1. {\rm\bigbreak\vskip-\parskip\noindent{\bf Remark #1.}\quad}

\def\Example #1. {\rm\bigbreak\vskip-\parskip\noindent{\bf Example #1.}\quad}

\def\Problems #1. {\bigbreak\vskip-\parskip\noindent{\bf Problems #1.}\quad}
\def\Problem #1. {\bigbreak\vskip-\parskip\noindent{\bf Problem #1.}\quad}

\def\Conjecture #1. {\bigbreak\vskip-\parskip\noindent{\bf Conjecture #1.}\quad}

\def\Proof#1.{\rm\par\ifdim\lastskip<\bigskipamount\removelastskip\fi\smallskip
            \noindent {\bf Proof.}\quad}

\def\Axiom #1. {\bigbreak\vskip-\parskip\noindent{\bf Axiom #1.}\quad\it}

\def\Satz #1. {\bigbreak\vskip-\parskip\noindent{\bf Satz #1.}\quad\it}

\def\Korollar #1. {\bbr\vskip-\parskip\nin{\bf Korollar #1.} \quad\it}

\def\Bemerkung #1. {\rm\bigbreak\vskip-\parskip\noindent{\bf Bemerkung #1.}
\quad}

\def\Beispiel #1. {\rm\bigbreak\vskip-\parskip\noindent{\bf Beispiel #1.}\quad}
\def\Aufgabe #1. {\rm\bigbreak\vskip-\parskip\noindent{\bf Aufgabe #1.}\quad}

\def\Beweis#1. {\rm\par\ifdim\lastskip<\bigskipamount\removelastskip\fi
           \smallskip\noindent {\bf Beweis.}\quad}

\nopagenumbers

\def\date{ }

\def\title{Title ??}
\def\author{Author ??}

\def\thanks#1{\footnote*{\eightrm#1}}

\def\rightheadline{\hfil{\eightrm\author}\hfil\tenbf\folio}
\def\leftheadline{\tenbf\folio
\hfil{\eightrm\title}\hfil}
\headline={\vbox{\line{\ifodd\pageno\rightheadline\else\leftheadline\fi}}}

\def\firstheadline{}
\def\firstfootline{\cen{\rm\folio}}

\def\seite #1 {\pageno #1
               \headline={\ifnum\pageno=#1 \firstheadline
               \else\ifodd\pageno\rightheadline\else\leftheadline\fi\fi}
               \footline={\ifnum\pageno=#1 \firstfootline\else{}\fi}}

%%%THIS IS THE MACRO LEFTSPACE.TEX %%%TO THD VIA WAFRUPP
\newdimen\dimenone
 \def\checkleftspace#1#2#3#4{%DIESER MACRO STAMMT VON APPELT
 \dimenone=\pagetotal%#1=Skip vorher,#2=Font,#3=Text,#4=Skip nachher  
 \advance\dimenone by -\pageshrink   %testen ob Titel noch mit Gewalt auf Seite 
                                                                          %geht
 \ifdim\dimenone>\pagegoal          %nacha tua nix-- gewoehnliche Outputroutine 
   \else\dimenone=\pagetotal
        \advance\dimenone by \pagestretch
        \ifdim\dimenone<\pagegoal
          \dimenone=\pagetotal
          \advance\dimenone by#1         %addieren Skip vor Ueberschrift (=#1)
          \setbox0=\vbox{#2\parskip=0pt                %#2 ist gewaehlter Font
                     \hyphenpenalty=10000
                     \rightskip=0pt plus 5em
                     \noindent#3 \vskip#4}    %#3=Ueberschrift,#4=skip nachher
        \advance\dimenone by\ht0
        \advance\dimenone by 3\baselineskip   
        \ifdim\dimenone>\pagegoal\vfill\eject\fi
          \else\eject\fi\fi}

%%% OUR HEADLINE MACROS LOOK LIKE THIS USING THIS MACRO

\def\sectionheadline #1{\bigbreak\vskip-\lastskip
      \checkleftspace{1.1cm}{\bf}{#1}{\bigskipamount}
         \vbox{\vskip1.1cm}\cen{\bfone #1}\bsk}

\def\lsectionheadline #1 #2{\bigbreak\vskip-\lastskip
      \checkleftspace{1.1cm}{\bf}{#1}{\bigskipamount}
         \vbox{\vskip1.1cm}\cen{\bfone #1}\msk \cen{\bfone #2}\bsk}

\def\lchapterheadline #1 #2{\bigbreak\vskip-\lastskip\indent\vskip3cm
                       \cen{\bftwo #1} \msk \cen{\bftwo #2} \gsk}
\def\llsectionheadline #1 #2 #3{\bigbreak\vskip-\lastskip\indent\vskip1.8cm
\cen{\bfone #1} \msk \cen{\bfone #2} \msk \cen{\bfone #3} \nobreak\bsk\nobreak}

%\def\[#1 #2\par{\hbox{\vtop{\hsize = 2.5 true cm \nin [#1]\hfill}
%\vtop{\hsize = 12.0 true cm \nin #2\penalty10000\llap.}}
%\vbox{\vskip.3\baselineskip}}  
% parameters for hsize are percentages !!! f.e. 0.2 + 0.8 = 1.0

\newtoks\literat
\def\[#1 #2\par{\literat={#2\unskip.}%
\hbox{\vtop{\hsize=.15\hsize\nin [#1]\hfill}
\vtop{\hsize=.82\hsize\nin\the\literat}}\par
\vskip.3\baselineskip}

\def\references{
\sectionheadline{\bf References}
\frenchspacing

\entries\par}

\mathchardef\emptyset="001F

\def\firstpage{\nin
{\obeylines \parindent 0pt }
\vskip2cm
\centerline{\bfone\title}
\gsk
\centerline{\bf\author}
\vskip1.5cm \rm}

% END OF LIEMACS.TEX

\def\date{jan 9, 2007}
\overfullrule = 0pt 

% \pageno=1

\def\title{Homotopes and Conformal Deformations of Symmetric spaces}

\def\author{Wolfgang Bertram}
% F\"ur endg\"ultige Version dieses Makro l\"oschen
%\def\leftheadline{\tenbf\folio\hfil habb/Contra.tex, \eightrm\date}

\def\Svar{\mathop{\rm Svar}\nolimits}
\def\Gras{\mathop{\rm Gras}\nolimits}

\def\YY{{\cal Y}}
\def\QQ{{\cal Q}}
\def\XX{{\cal X}}

\def\ZZ{{\cal Z}}

\def\U{{\cal U}}

\def\PPP{{\Bbb P}}

\def\v{{\frak v}}

\def\Gras{{\cal G}ras}

\def\PE{\mathop{\rm PE}\nolimits}

\def\Gl{\mathop{\rm GL}\nolimits}

% Strukturabbildungen
\def\S{{\bf \Sigma}}
\def\P{{\bf \Pi}}

\firstpage

\bigskip
\centerline{Wolfgang Bertram }
\centerline{
 Institut Elie Cartan (UMR 7502,  CNRS, INRIA)}
\centerline{
Universit\'e Henri Poincar\'e  -- Nancy-Universit\'e}
\centerline{B.P. 239}
\centerline{F - 54506 Vand\oe uvre-l\`es-Nancy Cedex} 
\centerline{e-mail: bertram@iecn.u-nancy.fr }

\vskip 10mm

\nin {\bf Abstract.} 
{\it Homotopy} is an important feature of associative and Jordan
algebraic structures: such structures always come in families whose
members need not be isomorphic among each other, but still share many
important properties.
One may regard homotopy as a special kind of {\it deformation}
of a given algebraic structure. In this work, we investigate the
global counterpart of this phenomenon on the geometric level of
the associated symmetric spaces -- on this level,
homotopy gives rise to 
{\it conformal deformations of symmetric spaces}.
These results are valid in arbitrary dimension and over general
base fields and -rings.

\bigskip
\nin {\bf Contents:}

\ssk
0. Introduction 

1. Homotopes of Lie- and Jordan algebraic structures

2. Deformations of quadratic prehomogeneous symmetric spaces

3. Deformations of generalized projective geometries

4. Deformations of polar geometries

5. Behavior of geometric structures under deformations

6. Some examples

Appendix A: Symmetric spaces

Appendix B: Generalized projective geometries

Appendix C: Globalization of structural maps

\vskip 10mm

\nin {\bf AMS-classification:}
17C37, % Associated geometries
32G99, % 32 Deformations of analytic structures :
      %  None of the above, but in this section
Secondary:
17C27, % Idempotents, Peirce decompositions
%17C30 Associated groups, automorphisms
%17C36 Associated manifolds
%17B60 Lie (super)algebras associated with other structures 
%         (associative, Jordan, etc.)
14D99 % 14Dxx Families, fibrations: None of the above, but in this section

\vskip 8mm

\nin {\bf Key words:} Homotope, isotope, Jordan algebras (-triple systems,
-pairs), Lie triple system, symmetric space, generalized projective
(and polar) geometries.

\vfill \eject

\sectionheadline
{0. Introduction}

\nin
{\bf 0.1. Conformal deformations.}
{\it Deformations} of Lie and other algebras have been much studied since
the fundamental work of Gerstenhaber [G64], and  related constructions on
the level of Lie groups (``contractions'') have attracted the interest of
physicists even earlier ([IW53]; see also [DR85]).
Following Gerstenhaber, one usually considers a given Lie algebra structure
$\mu$ on $\g$ as a point of the variety $\cal L$ of {\it all} Lie
algebra structures on $\g$, and one studies the structure of $\cal L$ 
on a neighborhood of $\mu$. In order to cover the case of base fields different
from $\R$, one considers {\it formal deformations} $\mu_t =
\sum_{m=0}^\infty t^m \mu_m$, and then relates conditions on the formal
derivatives of such one-parameter families of deformations to Lie algebra
cohomology.

\ssk
In this work we consider a special kind of deformations which we call
{\it conformal} in order to distinguish them from the general and more
formal approaches mentioned above. These deformations have the advantage
to lift always to the space level  and thus to define global
 deformations (or ``contractions'')  of geometric spaces,
given by explicit formulas which ensure a good control on the simultaneous
deformation of geometric properties of the spaces in question,
regardless of dimension and nature of the base field. For this
reason we hope that our approach may turn out to be useful for
harmonic analysis (in the spirit of [DR85]; cf.\ [FaPe04], [Pe02] for
topics to which our methods may apply) and maybe for physics (the 
contraction of de Sitter or anti-de Sitter models of general relativity
to Minkowski-space of special relativity is a special case of conformal
deformation in our sense, but many other cases, both finite and
infinite dimensional and potentially interesting for
physics, are also included in our setting).

\ssk
The possibility of conformal deformations of a geometric space
(symmetric space or Lie group) relies on the existence of an additional
structure which we call {\it generalized conformal structure}.
Algebraically, this means that we consider symmetric spaces or Lie groups
which lie in the  the image of the {\it Jordan-Lie functor} (see Section
1.6 for its definition):
there is an ``overlying'' Jordan algebraic structure from which one can
recover the Lie algebraic structure, and this overlying Jordan
structure corresponds
geometrically to a kind of conformal structure. 
Now, on the Jordan algebraic level there exists a notion of {\it homotopy} --
Jordan algebraic structures $J$ always come in natural
families $(J_i)_{i \in I}$, where the index set $I$ typically is a 
vector space or a certain algebraic variety; the Jordan structures $J_i$ are
then called {\it homotopes} of $J$, and if in turn $J$ can be obtained
as a homotope of $J_i$, then $J_i$ is called an {\it isotope} of $J$.
In our previous work [Be02], the relation between  geometric structures
belonging to isotopic Jordan structures
could be fully understood, but it was not clear how to understand
``deformations'' or ``contractions'' to more singular values, where
$J_i$ is only homotopic, but 
 no longer isotopic to $J$. In the present work we close this gap.

\ssk
To our knowledge, there is no interesting notion of homotopy in a purely
Lie-theoretic setting, and therefore it seems quite hard to recover the
present results without using Jordan-theory (if one wishes to avoid
case-by-case calculations). On the other hand, once the geometric framework
is clarified, some results given in this paper could be used
to derive in a geometric way many of the identities and properties of
Jordan pairs from [Lo75]. We illustrate this remark by giving independent
proofs of some of our results, one proof based on the algebraic
results from [Lo75], and
another based 
 on the geometric framework of {\it generalized projective geometries}
([Be02]; cf.\	Appendix B).

\msk
\nin
{\bf 0.2. Example: associative algebras and general linear groups.}
For readers who  are not experts in Jordan theory, let us give
a very simple
example that is well-suited to illustrate
what we mean by ``conformal deformation'', and by   ``homotopy'':
if $G$ is any group with unit element $\1$, then we can define new group
laws on $G$ by choosing an arbitrary element $a \in G$ and defining
$$
x \cdot_a y := x a y.
$$
The unit element of $(G, \cdot_a)$ is then $a^{-1}$.
Of course, all these group laws are isomorphic among each other:
left translation $l_a:G \to G$ is an isomorphism from $(G,\cdot_a)$
onto $(G,\cdot)$ (one may also use right translation), and hence
this construction seems not to be very interesting.
However, things would start to become more interesting, if one could
take a ``limit'' for $a$ converging to some element in the ``boundary
of $G$''; then $(G,\cdot_a)$ should converge to some 
deformed group  which is no longer isomorphic to $G$ (a  ``homotope''
of $G$).

In order to make this idea more precise, we need some additional structure.
For instance, assume that $G = A^\times$ is the unit group of an
associative unital algebra $A$ over a commutative  ring  $\K$.
As explained above, $G:=A^\times$ is a group with product $\cdot_a$,
for any $a \in A^\times$. We want to allow $a$ to ``converge'' to some non-invertible
element of $A$ -- the problem being that, although the {\it homotope algebra}
$(A,\cdot_a)$ with 
product $x \cdot_a y = xay$ is an
associative algebra for all $a \in A$, it will not have a unit element
as soon as $a$ is non-invertible.
In order to circumvent this problem, 
  we translate, for invertible $a$, the unit element
$a^{-1}$ by the translation $\tau=\tau_{-a^{-1}}$ to the origin $0 \in A$, 
that is, for invertible $a$ we consider the group law
$$
\eqalign{
x \diamond_a y&
= \tau (\tau^{-1}(x) \cdot_a \tau^{-1} (y)) \cr
& = (x+a^{-1}) a (y + a^{-1}) - a^{-1} \cr
& = x a y + y + x \cr}
\eqno (0.1)
$$
which is defined on the set 
$$
\eqalign{
G_a & =  \{ x \in A | \, a^{-1} + x \in A^\times \} \cr
& = \{ x \in A | \, \1 + a x \in A^\times \}. \cr}
\eqno (0.2)
$$
Now it is easy to take 
the ``limit as $a$ becomes non-invertible'': 
one checks that the last expressions of
(0.1) and (0.2) define a group $(G_a,\diamond_a)$ for
any element $a \in A$.
(The inverse of $x \in G_a$ is
$y= - x (\1+ax)^{-1}$, as is seen by direct calculation.)
 For $a \in A^\times$ these groups are, by construction,
all isomorphic to $G$, but for non-invertible $a \in A$ they are not:
they are ``homotopes'' , but not isotopes of $G$.
% If $A$ is a topological algebra over a 
%topological ring, then every element $a$ belonging to the closure of
%$A^\times$ will give rise to a deformation of the group $A^\times$.
In particular, for $a=0$, we get the translation group $(A,+)$ as a deformation
of $A^\times$.
Also, it is clear that, if $a,b \in A$ belong to the same
$A^\times$-orbit (under left, or under right action), then $G_a$ and $G_b$
are isomorphic. Therefore, if we can classify the orbits of the
action of $A^\times \times A^\times$ on $A$, then we can also classify
the various groups $G_a$ obtained by the preceding construction.
For instance, this is possible in the case where
$A = M(n,n;\K)$ is the algebra of square matrices
over a field $\K$, and then a complete description of the groups $G_a$ 
is not difficult (see Section 5.1). 
%
% remark: ``flocks'' , cf. Pruefer, Baer, Certaine, book of Sabinin.
%
Finally, one remarks that,
on the level of Lie algebras, the corresponding situation is very simple:
whenever the notion of a ``Lie algebra of $G_a$" makes sense, then 
this is simply  $A$ with the $a$-deformed Lie bracket
$$
[x,y]_a = xay - yax. \eqno (0.3)
$$

\msk \nin {\bf 0.3. Plan of this work.}
In the first chapter we recall the definition of the three basic
Jordan theoretic categories and the notions of homotopy we are interested in:

\ssk
\item{UJA} Unital (linear) Jordan algebras $(J,\bullet)$
 (Section 1.7), with the $a$-homotopes $\bullet_a$, $a \in J$.
\item{JP} (Linear) Jordan pairs
 $(V^+,V^-)$ (Section 1.4), with the $a$-homotopes $V^+_a$, $a \in V^-$. 
\item{JTS} (Linear) Jordan triple systems $(V,T)$, with the $\alpha$-homotopes
$(V,T_\alpha)$, where $\alpha$ is a {\it structural transformation}
(Lemma 1.13).

\ssk
\nin In Chapter 2 we describe the Jordan theoretic analog of
 the associative construction given above (Prop.\ 2.2,
related to the $a$-homotope $\bullet_a$, and Th.\ 2.3, related to
$V^+_a$). In Chapter 3 we give an independent and geometric version
of these results (Th.\ 3.1 and Th.\ 3.3).
Finally, in Chapter 4 we construct the geometric version of deformations
corresponding to the $\alpha$-homotope $T_\alpha$ (Th.\ 4.3).
In some sense, this is the most general (but also most abstract) result,
since one may choose $\alpha = Q(a)$, the quadratic operator corresponding
to $a \in V^-$, in order to recover the $a$-homotope of $T$.
In view of applications in harmonic analysis or physics (see above), we
describe in Chapter 5 the behavior of geometric structures (for instance,
affine connections or duals of metrics) under
conformal deformations -- for brevety, whenever calculations are essentially
the same as in
[Be00] and [Be02], we refer to these works for further details. 
In Chapter 6 we make our abstract results more explicit for some cases,
namely for $a$-deformations belonging to {\it idempotents} or
{\it tripotents} (in the semisimple
case of finite dimension over a field this makes
a classification possible), and for the case of {\it Grassmannians}.
In the appendices,
we collect some basic facts on reflection spaces and
symmetric spaces (A), on generalized projective geometries (B) and
on structural transformations (C); the simple Lemma C.3 on globalizations
of such transformations is a main ingredient to be used in the proof of
Theorem 4.3.

\msk \nin
{\bf Acknowledgement.}  I thank Ottmar Loos for helpful comments.
%In achieving this work I benefitted of a teaching
%release which was possible thanks to a common effort of students
%und university staff of the university Henri Poincar\'e 
%and which a gratefully acknowledge.

\msk
\nin {\bf Notation.}
Throughout, we work over a commutative base ring $\K$
in which $2$ and $3$ are invertible and 
hence use the concepts {\it linear} Jordan pairs, algebras and triple systems.

\sectionheadline
{1. Homotopes of Lie and Jordan algebraic structures}

\msk \nin
{\bf 1.1. Scalar homotopes and c-duality.}
If $V$ is an algebra over $\K$ (Lie or other), with product $m(x,y)$,
then for any scalar $\lambda \in \K$, we may define a new product
$m_\lambda(x,y) := \lambda m(x,y)$, called a {\it scalar homotope of$(V,m)$}. 
Then $\phi:=\lambda \id_V$ is a homomorphism from $(V,m_\lambda)$ to 
$(V,m)$.
 In particular, the algebras $(V,m_\lambda)$
for invertible $\lambda$ are all isomorphic
among each other. Thus, if $\K$ is a field, this construction
gives us nothing interesting. 

Now  let $(V,R)$ be a triple system (i.e., $R:V^3 \to V$ is a trilinear map).
Then
$R_\lambda(X,Y,Z):=\lambda R(X,Y,Z)$
defines a new triple product on $V$, again called a {\it scalar homotope}. 
But now $\lambda \id_V$ is a homomorphism from $(V,R_{\lambda^2})$ to
$(V,R)$, and hence we can only conclude that
all homotopes $(V,R_\mu)$, with $\mu$ a square of an invertible
element in $\K$, are isomorphic to $(V,R)$.
For instance, if $\K=\R$, we get a new  triple system
$(V,R_{-1})$ which is called the {\it $c$-dual of $R$}.
For a general field, we say that the $(V,R_\lambda)$ with invertible 
$\lambda$ are {\it scalar isotopes
of $(V,R)$}; they need not be isomorphic to $(V,R)$. 

The aim of this chapter is to recall some basic results showing that
for Jordan algebraic structures there are interesting generalizations
of the notion of scalar homotopy.

\msk \nin
{\bf 1.2. Graded Lie algebras.}
Let $(\Gamma,+)$ be an abelian group. A {\it $\Gamma$-graded
Lie algebra} is a Lie algebra $\g$ over $\K$ with a decomposition
$\g = \bigoplus_{n \in \Gamma}
\g_n$ such that  $[\g_m,\g_n] \subset \g_{n+m}$ for all $n,m \in \Gamma$.
The $\Z/(2)$-graded Lie algebras correspond bijectively to {\it symmetric
Lie algebras}, i.e., to Lie algebras together with an automorphism of order
$2$. 
A  {\it $3$-graded} Lie algebra is a  $\Z$-graded Lie algebra such that
$\g_n=0$ if $n \not= -1,0,1$.
Note that then $\g = \g_0 \oplus (\g_1 \oplus \g_{-1})$ gives rise
to a $\Z/(2)$-grading, so $\g$ is then also a symmetric Lie algebra (of a 
rather special type).

\msk
\nin
{\bf 1.3. Lie triple systems.} A Lie triple system (LTS) is a $\K$-module $\q$
together with a trilinear map 
$$
\q \times \q \times \q \to \q, \quad (X,Y,Z) \mapsto [X,Y,Z]
$$
such that (where we use also the notation $R(X,Y,Z):=R(X,Y)Z:=-[X,Y,Z]$,
alluding to the interpretation of this expression as a {\it curvature tensor}
in differential geometry)

\ssk
\item{(LT1)} $R(X,X)=0$
\item{(LT2)} $R(X,Y)Z+R(Y,Z)X+R(Z,Z)Y =0$ (the Jacobi identity)
\item{(LT3)} For $D:=R(X,Y)$ we have 
$D R(U,V,W)=R(DU,V,W)+R(U,DV,W)+R(U,V,DW)$.

\ssk
\nin If $\g$ is a $\Z/(2)$-graded Lie algebra, then $\g_{-1}$ 
with $[X,Y,Z]:=[[X,Y],Z]$ becomes an LTS, and every LTS arises in this way
since one may reconstruct a Lie algebra from an LTS via the {\it
standard imbedding} (due to Lister, cf.\ [Lo69]).

\msk \nin
{\bf 1.4. Linear Jordan pairs.}
If $\g = \g_1 \oplus \g_0 \oplus \g_{-1}$ is  a $3$-graded Lie algebra,
 then we let $V^\pm := \g_{\pm 1}$ and define trilinear maps by
$$
T^\pm:  V^\pm \times V^\mp \times V^\pm \to V^\pm, \quad
(x,y,z) \mapsto  [[x,y],z].
$$
The maps $T^\pm$ satisfy the following identies:

\ssk
\item{(LJP1)} $T^\pm(x,y,z)=T^\pm(z,y,x)$
\item{(LJP2)} 
$T^\pm(u,v,T^\pm(x,y,z))=
T^\pm (T^\pm(u,v,x),y,z) -
T^\pm(x,T^\mp(v,u,y),z) +
T^\pm(x,y,T^\pm(u,v,z))$

\ssk
\nin In general, a pair of $\K$-modules $(V_+,V_-)$ with trilinear maps
$T^\pm: V_\pm \times V_\mp \times V_\pm \to V_\pm$ is called a 
{\it linear Jordan pair} if (LJP1) and (LJP2) hold.
Every linear Jordan pair arises by the construction just described.
%The reconstruction of the Lie algebra $\g$ from a Jordan pair $(V^+,V^-)$
% is essentially the standard imbedding in disguise since
%$\g_1 \oplus \g_{-1}$ is a Lie triple system.

\msk \nin
{\bf 1.5. Jordan triple systems.}
An {\it involution of a $3$-graded Lie algebra} is an automorphism $\theta$
of order
$2$ such that $\theta(\g_\pm)=\g_\mp$ and $\theta(\g_0)=\g_0$.
Then we let $V:=\g_+$ and define
$$
T:V \times V \times V \to V, \quad (X,Y,Z) \mapsto [[X,\theta Y],Z].
$$
Then $T$ satisfies the identities (LJT1) and (LJT2) obtained from
(LJP1) and (LJP2) by omitting the indices $\pm$. A $\K$-module
together with a trilinear map satisfying these identities is called a
 {\it (linear) Jordan triple system} (JTS). Every linear Jordan triple system
arises by the construction just described.
One may also define the notion of an {\it involution of a Jordan pair};
then Jordan triple systems are the same as Jordan pairs with involution
(see [Lo75]). For any JTS $(V,T)$, the Jordan pair $(V^+,V^-)=(V,V)$ with
$T^\pm=T$ is called the {\it underlying Jordan pair}. 

\msk \nin
{\bf 1.6. The Jordan-Lie functor.}
If $T$ is a JTS on $V$, then
$$
[X,Y,Z ]= T(X,Y,Z) - T(Y,X,Z)
$$
defines a LTS on $V$. This defines a functor from the category of JTS to
the category of LTS over $\K$, which we call the {\it Jordan-Lie functor}
 (cf.\ [Be00]). 
The LTS $(V,R)$ belongs the $\theta$-fixed subalgebra $\g^\theta$ of $\g$ 
whose involution is
induced by the canonical involution of the $3$-graded Lie algebra $\g$.

\msk \nin
{\bf 1.7. Jordan algebras.}
A {\it Jordan algebra} is a $\K$-module $J$ with a bilinear and 
commutative product
$x \bullet y$ such that the identity

\item{(J2)} $x \bullet (x^2 \bullet y) = x^2 \bullet (x \bullet y)$

\nin holds (cf.\ [McC04]; 
the reader not familiar with Jordan theory may throughout think
of the example $J=A^\tau$, the set of fixed points of 
an anti-automorphism $\tau$ of an associative algebra $A$, with 
bilinear product
$x\bullet y={x  y + y  x \over 2}$.) 
Jordan algebras always come in families:
for any Jordan algebra $J$ and element $a \in J$, the {\it homotope}
$J_a$ is defined as follows. Writing $L(x)$ for the
 left multiplication operator
$L(x) y =x\bullet y$,  we can define a new  product on $J$ by
$$
x \bullet_a y:= (L(xa) + [L(x),L(a)])y = (x \bullet a) \bullet y + x
 \bullet (a \bullet y) - a \bullet (x \bullet y).
$$
One proves that $J_a=(J,\bullet_a)$ is a again a Jordan algebra
 (see [McC04]), called the {\it $a$-homotope of $(V,\bullet)$}, 
 and that $T(x,a,y) :=x \bullet_a y$ defines a Jordan
 triple product on $J$
(see [Be00, Ch.\ II] for a geometric proof of this fact). Thus $(J,J)$
is a Jordan pair, but not all Jordan pairs are of this form.

\msk \nin {\bf 1.8. Jordan pairs with invertible elements.}
For any Jordan pair $(V^+,V^-)$ and every $a \in V^-$, $V^+$ with the
bilinear product 
$$
x \bullet_a y:= {1 \over 2} T^+(x,a,y)
$$
is a Jordan algebra, again called the {\it $a$-homotope algebra}.
This algebra is unital if, and only if, the element $a$ is 
{\it invertible in $(V^+,V^-)$}, which means that
 the  operator
$$
Q(a):=Q^-(a):={1 \over 2} T^-(a,\cdot,a): V^+ \to V^-
$$
is invertible; then $a^{-1}:=(Q(a))^{-1}(a)$ is the unit element of $J_a$.
Every unital Jordan algebra arises in this way from a Jordan pair
having invertible elements, after a choice of such an element.

\msk \nin {\bf 1.9. Invertibility.}
In a Jordan algebra $(J,\bullet)$, the quadratic operator $U_x$ is defined by 
$$
U_x y  =  (2 L_x^2 - L_{x^2}) y = 2 x \bullet (x \bullet y) - (x \bullet
x) \bullet y.
$$
An element $x \in J$ is called {\it invertible} if $U_x$ is invertible;
let
$$
J^\times = \{ x \in J | \, U_x \in \Gl(J) \},
$$
be the set of invertible elements.
For such $x$, we define the {\it inverse} $x^{-1}=U_x^{-1}x$ and call the map 
$j:J^\times  \to J$, $x \mapsto x^{-1}$ the {\it Jordan inverse}.
Then the following identities hold (cf.\ [McC04, p.\ 200 and
201]):

\ssk
\item{(SB1)} $U_{U_x y} = U_x U_y U_x$ 
(the ``fundamental formula'')
\item{(SB2)} $U_{x^{-1}} = (U_x)^{-1}$
\item{(SB3)} $x^{-1} = U_x^{-1}x$

% In Section 2.1 we give a geometric interpretation of these identities.

\msk \nin {\bf 1.10. Quasi-invertibility.}
In general, a Jordan pair $(V^+,V^-)$ does not have invertible elements.
The important notion is the one of {\it quasi-inverse}: let $$
\eqalign{
B(x,y) :&= \id_V - T(x,y) + Q(x) Q(y) \cr
\tilde \tau_y (x) :&= x^y := B(x,y)^{-1} (x-Q(x)y) \cr}
$$
 the {\it Bergman operator}, resp.\ the
{\it quasi-inverse} (the latter provided that $(x,y)$ is 
{\it quasi-invertible}, which means that $B(x,y)$ is invertible).
If $x$ is invertible and $(x,y)$ quasi-invertible, then we have 
$\tilde \tau_y(x)=(x^{-1}+y)^{-1}$. 

%
%In case $(V^+,V^-)=(J,J)$ comes from a unital Jordan algebra $V$,
%we can also write
%$\tilde \tau_y = j \circ \tau_{-y} \circ j$ where $j$ is the Jordan inverse.
%For $y=1$ (the unit element), the identities
%$\tilde \tau_1 = j \circ \tau_{-1} \circ j$ 
%and
%$\tau_{-1} \circ j \circ \tau_1 =  \tilde \tau_1 \circ (-\id)$, 
%where $j$ is the Jordan inverse, can be interpreted via 
%an action of the group $\Sl(2,\K)$ on $J$ by fractional linear maps via
%$$
%\pmatrix{a & b \cr c & d \cr}.x = (ax+b 1)(cx + d1)^{-1}.
%$$
%Then the maps $\tau_1$, $\tilde \tau_1$, $j$ and $-\id$ correspond to the
%matrices
%$$
%\pmatrix{1 & 1 \cr 0 & 1 \cr}, \quad
%\pmatrix{1 & 0 \cr -1 & 1 \cr}, \quad
%\pmatrix{0 &1 \cr 1 & 0 \cr}, \quad
%\pmatrix{-1 & 0 \cr 0 & 1 \cr},
%$$
%and one checks the corresponding equations by a direct matrix computation.
%

\msk \nin
{\bf 1.11. Structural transformations
 of Jordan triple systems and associated homotopes.}
For Jordan triple systems, there is a much richer supply of
homotopes than for Lie triple systems: let $(V,T)$ be a
 Jordan triple system.
A {\it structural transformation of $(V,T)$} is
an endomorphism  $\alpha$ of $V$ with the property
$$
\forall x,y,z \in V: \quad T(\alpha x, y, \alpha z)=\alpha
T(x,\alpha y,z). 
\eqno (1.1)
$$

\Lemma 1.12. If $\alpha$ a structural transformation, then the formula
$$
T^{(\alpha)}(x,y,z):= T(x,\alpha y,z)
\eqno (1.2)
$$
defines a Jordan triple product $T_\alpha$ on $V$, called
the {\rm $\alpha$-homotope of $T$}. The quadratic map and the Bergmann operator
of $T_\alpha$ are given by
$$
Q_\alpha(x)=Q(x) \circ \alpha, \quad \quad
B_\alpha(x,y)=B(x,\alpha y).
$$

\Proof.  It is clear that $T_\alpha$ satisfies (LJT1). For the proof
of (LJT2) take the
identity (LJT2) for $T$ and replace the middle elements $v$ and $y$ by
$\alpha v$ and $\alpha y$ and get
$$
T(u,\alpha v,T(x,\alpha y,z))=
T (T(u,\alpha v,x),\alpha y,z) -
T^\pm(x,T(\alpha v, u,\alpha y),z) +
T^\pm(x,\alpha y,T^\pm(u,\alpha v,z))
$$
Then apply the defining relation of $\alpha$ to the second term on the
right hand side, and one gets the identity (LJT2) for $T_\alpha$. 
The quadratic operator for this JTS clearly is
$Q_\alpha(x)=Q(x) \circ \alpha$, and it follows that
$$
\eqalign{
B_\alpha(x,y)& =\id-T_\alpha(x,y)+Q_\alpha(x)Q_\alpha(y)=
\id - T(x,\alpha y)+Q(x) \circ \alpha Q(y) \circ \alpha \cr
& =\id - T(x,\alpha y) + Q(x) Q(\alpha y) =
B(x,\alpha y).   \cr}
\qeddis

\nin The set of structural transformations of $(V,T)$,
$$
\Svar(T):=\{ \alpha \in \End(V) | \,
\forall x,y,z \in V: \quad T(\alpha x, y, \alpha z)=\alpha
T(x,\alpha y,z) \}
$$
is called the {\it structure variety of the Jordan triple system $T$}.
Note that, by the fundamental formula $Q(Q(x)z)=Q(x)Q(z)Q(x)$,
all quadratic maps $Q(x)$ belong to $\Svar(T)$, and so do the Bergmann
operators of the form $B(x,x)$ for all $x \in V$ (by Bergmann structurality
$Q(B(x,y)z)=B(x,y)Q(z)B(y,x)$, cf.\ [Lo75, JP26]).
Also, $\Svar(T)$ contains all involutive automorphisms of $T$,
together with their negatives. More generally, $\Svar(T)$ is a cone, i.e.,
invariant under taking scalar multiples.

\sectionheadline
{2. Deformations of quadratic prehomogeneous symmetric spaces}

\nin {\bf 2.1. The symmetric space of invertible elements in a Jordan
algebra.}
Assume $J$ is a Jordan algebra over a ring $\K$ with Jordan product
$x \bullet y$ and unit element $1$. 
The set $J^\times$ of invertible
elements of $J$ is stable under the maps
$j(x)=x^{-1}$ (the Jordan inverse) and $U_x$ for every $x \in J^\times$.
Let $\QQ(x):=U_x$. Then the identities 
(SB1) -- (SB3) from Section 1.9 are exactly the identities
(SB1) -- (SB3) from Appendix A.1, and therefore
Prop.\ A.2 implies that $J^\times$ carries a natural
structure of a reflection space given by
$$
\mu(x,y)=\sigma_x(y)=U_x (y^{-1}) = U_x U_y^{-1}y.
$$
The differential at $1$ of the symmetry $\sigma_1 = j$ can be calculated
algebraically by using scalar extension by dual numbers $\K[\eps]=
\K \oplus \eps \K$, $\eps^2=0$. It is easily checked that, in
$J \otimes_\K \K[\eps]=J \oplus \eps J$, we have
$(\1 + \eps v)^{-1} = \1 - \eps v$, and it follows that the tangent map
of $j$ et the identity is $- \id_J$. 
If $a \in J^\times$ is another base point, the maps $\QQ$ and $j$
are replaced by the quadratic map, resp.\ Jordan inverse in the
homotope algebra $(J,\bullet_{a^{-1}})$ (which has $a$ as unit element), 
and hence we have also $T_a \sigma_a = - \id_J$. 
Thus also property (S4) defining a symmetric
space is satisfied (in an algebraic sense), and hence $J^\times$ is a
symmetric space. 
We call such spaces {\it quadratic pre-homogenoues symmetric spaces}
(cf.\ [Be00, Chapter II]).
Now we describe what happens when we translate the space $J^\times$ such
that a base point $a$ is translated into the origin $0 \in J$.
Recall from Section 1.10 the definition of the Bergmann-operator and the
quasi-inverse. 

\Proposition 2.2.
We equip $J^\times$ with its symmetric space structure
$\mu(x,y)=\sigma_x(y)=U_x y^{-1} $, let $a \in J^\times$ and
$\tau = \tau_{-a^{-1}}:J \to J$ be translation by $-a^{-1}$.
Then the translated symmetric space $\tau (J^\times)$ is given by
$$
\tau(J^\times)=\{ x \in J \vert \, B(x,-a) \in \Gl(J) \}
$$
with product map
$$
\eqalign{
\tau \bigl( \mu (\tau^{-1}(x),  \tau^{-1}(y) \bigr) 
& = 2x+U_x a+B(x,-a) \tilde \tau_{a}(-y) \cr} .
$$

\Proof. 
First of all, $x \in \tau (J^\times)$ if and only if
$a^{-1} +x$ is invertible, if and only if 
$U_{a^{-1} +x}$ is invertible; but since $U_a$ is assumed to be invertible,
we have the identity 
$$
B(x,a)= U_{x-a^{-1}} U_a \eqno (2.1)
$$ 
([Lo75, I.2.12])
and hence $x$ belongs to $\tau(J^\times)$ if and only if $B(x,-a)$ is
invertible, i.e., iff $(x,-a)$ is a quasi-invertible pair.

Now let us calculate the translated product map, or, which amounts to the 
same, via Prop.\ A.2, the new
quadratic map and inversion with respect to a base point.
With respect to the old base point $a^{-1} \in J^\times$, the quadratic map
was given by $\QQ(x)y=U_x U_a y$, and the symmetry was the Jordan inverse
$j_a$ in the Jordan algebra $J_a$ with unit element $a^{-1}$.
We have to show that transformation by $\tau$ leads to the new quadratic
map and new inverse
$$
\eqalign{
\tau \Big( \QQ \big(\tau^{-1}(x)\big) \tau^{-1}(y) \Big) 
& = 2x+U_x a+B(x,-a) y,        \cr
\tau ( \sigma_a (\tau^{-1}(x))   & = \tilde \tau_{a}(-x) . \cr} 
\eqno (2.2)
$$
The first equality in (2.2) follows from (2.1): 
$$
\eqalign{
U_{x+a^{-1}}  U_a  (y+a^{-1})  - a^{-1} 
& = B(x,-a) (y + a^{-1})  \cr
& = B(x,-a)y + a^{-1} + T(x,a,a^{-1}) + 
 U_x U_a a^{-1} \cr
& =  B(x,-a)y + a^{-1} + 2x +   U_x a   \cr}
$$
For the proof of the second equality from (2.2), note that the inverse
in $J_a$ is $j_a(x)=U_a^{-1} x^{-1}$, and hence
$$
\tau \circ j_a \circ \tau^{-1}(x)= -a^{-1} + U_a^{-1}(x+a^{-1})^{-1} 
= - a^{-1} + U_a^{-1}(a^{-x}) = (-x)^a
$$
by the symmetry principle for the quasi-inverse, see [Lo75, 3.3].
\qed

\msk
\nin 
Next we will show that the formul\ae{ } of the preceding proposition
define a symmetric space for any $a \in J$, invertible or not.
Note that for $a=0$ this is trivial: in this case we get the flat
symmetric space structure on $J$. 
The proof becomes more transparent in the more general framework
of {\it Jordan pairs} -- in this context, the following equivalent version
of the formulae from Prop.\ 2.2 holds:

\Theorem 2.3.
Assume $(V^+,V^-)$ is a linear Jordan pair over $\K$ and fix an element
$a \in V^-$. Let
$\U_a:= \{ x \in V^+ \vert \, B(x,-a) \in \Gl(V^+) \}$
be the set of elements $x \in V^+$ such that $(x,-a)$ is quasi-invertible.
Then, for all $x \in \U_a$, the set $\U_a$ is stable under the transformations
$$
\eqalign{
\QQ(x) &:= \tau_{2x+Q(x)a} \circ B(x,-a) ,\cr
\sigma_0 &:= - \tilde \tau_{-a} = \tilde \tau_a \circ (-\id), \cr}
$$
and $(\U_a,\QQ,\sigma_0)$ is a reflection space with base point $0$.
Moreover, the algebraic differential of the symmetry $\sigma_x$ at $x$,
defined by $\sigma_x(x + \eps v)=\sigma_x(x)+\eps D\sigma_x(x)v$,
equals $-\id_{V^+}$, and thus $\U_a$ is a symmetric space.

\Proof.
We use the homotope $J=V_a^+$ with product
$ x \bullet_a y = {1 \over 2} T(x,a,y)$
and adjoin a unit element $\hat J = \K \1 \oplus J$
(if $J$ happens to be unital, we forget the unit element of $J$).
For $x,y \in J$, we calculate in $\hat J$, 
$$
\eqalign{
U_x y & = (2 L_x^2 - L_{x^2})y = 2 x \bullet_a (x \bullet_a  y) - 
(x \bullet_a x)  \bullet_a y = Q^+(x)Q^-(a)y\cr
U_{\1 \oplus x}(y) & = 
y + 2 x  \bullet_a y + U_x y = (\id + T(x,a) + Q(x)Q(a))y =
B(x,-a)y,
\cr
U_{\1\oplus x}(\1\oplus y) & = U_{\1 \oplus x}(\1) + U_{\1\oplus x}(y) = 
(\1 + x)^2 + U_{\1\oplus x}y
 = \1 \oplus \big(2x +  x^2 + U_{\1\oplus x}y\big)\cr
& = \1 \oplus \big(2x + Q(x)a + B(x,-a)y \big).
\cr}
\eqno (2.3)
$$
The second equality shows that
$B(x,-a)$ is the restriction of $U_{\1\oplus x}$ to $J$, and hence
$(x,-a)$ is quasi-invertible in $(V^+,V^-)$ 
if and only if $\1+x$ is invertible in $\hat J$, i.e.,
$$
\1 + \U_a = \hat J^\times \cap (\1 + J)  .
$$
This set
is a symmetric subspace of $\hat V^\times$. 
In fact, the third equation of (2.3) shows that it is stable under the
quadratic map $\QQ$ (with respect to the origin $1$), and it is also
stable under inversion: if $(x,-a)$ is quasi-invertible, then
the inverse of $\1+x$ in $\hat J$ is given by (see [Lo75, I.3])
$$
(\1+x)^{-1} = \1 - \tilde \tau_{-a} (x) .
$$
Transferring the symmetric space structure from
$\1+\U_a$ to $\U_a$ via the map $\1 \oplus x \mapsto x$, 
we obtain the claimed formula.

Algebraically, the differential $D\sigma_0(0)$ of $\sigma_0$ at the origin may
be defined by scalar extension with dual numbers $\K[\eps]=
\K \oplus \eps \K$, $\eps^2=0$:
$$
\sigma_0(0+\eps v) = - B(\eps x,-a)^{-1} (\eps x + Q(\eps x)a)=
- (\id - \eps T(x,a)) \eps x =  - \eps x,
$$
whence $D\sigma_0(0)=-\id_{V^+}$.
Since $\sigma_x = \QQ(x) \circ \sigma_0$, it follows that
$D\sigma_x(x)=-\id_{V^+}$.
\qed

\msk \nin
Remark. By the preceding calculations, the fundamental formula 
$U_{\1-x}U_{\1-y}U_{\1-x}=U_{U_{\1-x}(\1-y)}$
in $\1-\U_a$ translates to the identity
$$
B(x,a)B(y,a)B(x,a) = B\bigl( 2x-Q(x)a+B(x,a)y,a \bigr),
$$
which in turn can be interpreted by saying that, for every
$a \in V^-$,  the map
$x \mapsto B(x,a)$ is a homomorphism of the reflection
space $\U_a$ into the automorphism group of $(V^+,V^-)$, 
seen as a reflection space.

\msk \nin {\bf 2.4. The topological case.}
Assume $\K$ is a topological ring with dense unit group.
A {\it topological Jordan pair} is a Jordan pair together with topologies
on $V^+$ and $V^-$ such that $V^\pm$ are
Hausdorff topological
$\K$-modules $V^\pm$ and the trilinear maps $T^\pm$ are continuous.
Following [BeNe05], we say that {\it $V^\pm$ 
satisfies the condition} (C2) if the following holds:

\ssk
\item{(C2)} for all $a \in V^\mp$, 
the sets $\U_a$ are open and the maps $\U_a \times V^\pm \to V^\pm$,
$(x,z) \mapsto B(x,-a)^{-1} z$ are continuous.

\Theorem 2.5.
Assume $(V^+,V^-)$ is a topological Jordan pair 
satisfying the condition {\rm (C2)}. 
\item{(1)}
The symmetric space $(\U_a,\mu)$ from {\rm Theorem 2.3}
is smooth, i.e., it is a symmetric space in the category of manifolds over $\K$.
\item{(2)} 
With respect to the base point $0$, the Lie triple system of the
symmetric space $(\U_a,\QQ,\sigma_0)$ is given by antisymmetrizing
the Jordan triple product 
$S(u,v,w):= T^+(u,Q(a)v,w)$ in the first two variables:
$$
[x,y,z]=S(x,y,z)-S(y,x,z).
$$

\Proof. (1) Smoothness of $\mu$ follows from
 [BeNe05, Prop.\ 5.2], and the Property (S4) of a symmetric space
 follows from the calculation of $D\sigma_0(0)$ given above
 (the differential in the sense of differential calculus is the same
 as the algebraic differential, cf.\ [Be06]).
 
(2) 
Once again, the easiest proof is by reduction to unital Jordan algebras: 
for a unital Jordan algebra $(J,e)$,  the Lie triple system of the
pointed symmetric space $(J^\times,e)$ is given by
$$
[x,y,z] = T(x,y,z)-T(y,x,z)  = 
[L_x,L_y] z = x \bullet (y\bullet z) - y\bullet ( x\bullet z)
\eqno (2.4)
$$
(see [BeNe05, Theorem 3.4] and its proof which is algebraic in nature).
For invertible $a$, the Lie triple system of $(J^\times,a)$ 
is then given by the same formula, but replacing $\bullet$ by $\bullet_a$
and $T$ by $T_a(x,y,z)=T(x,Q(a)y,z)=S(x,y,z)$. 

Starting now with a Jordan pair $(V^+,V^-)$ and arbitrary $a \in V^-$,
we apply the preceding remarks to
the unital Jordan algebra
$\hat J = \K \1 \oplus J$ with $J=V_a^+$, as in the proof of Theorem 2.3.
Then the Lie triple system of the pointed symmetric space $(\1+\U_a,\1)=
\hat J^\times \cap (\1 + J)$ is
given by restricting the triple product (2.4) (with $\bullet =
\bullet_a$) to $J=V_a^+$.
\qed

\nin The preceding result applies to all finite-dimensional and to all
Banach Jordan pairs over $\R$ or $\C$, but also to many other 
not so classical cases.
Note that the (C2)-condition assures in particular
that the sets $\U_a$ are not reduced to $\{ 0 \}$ (if $V^+ \not= \{ 0 \}$). 
As to Part (2) of the theorem,
it is a general fact that antisymmetrizing a Jordan triple
system in the first two variables gives a Lie triple system (the 
``Jordan-Lie functor", see Section 1.6).
The preceding construction works also in other contexts where 
some sort of differential calculus can be implemented (cf.\
Remark 10.1 in [BeNe05]), and it
can also be applied if we start with any non-unital Jordan algebra
instead of $J_a$; then $\U_a$ corresponds to the set of quasi-invertible
elements in this algebra.

%\msk Drawback of adjoining a unit :
%it works for fixed $a$; the deformation aspect gets lost.
%But it facilitates greatly the algebraic proof.
%The following geometric arguments in some sense give the proof if one
%does not use the adjoining of a unit.

\sectionheadline
{3. Deformations of generalized projective geometries}
                                                                   
\nin The geometric or  ``integrated" version of a Jordan pair $(V^+,V^-)$
is a {\it generalized projective geometry} $(\XX^+,\XX^-)$, see
Appendix B. We give an analog of Theorem 2.3 in this context; it 
generalizes [Be03, Section 6], where
only the case of Jordan algebras has been treated (cf.\ last part of the claim).

\Theorem 3.1.
Let $(\XX^+,\XX^-)$ be the generalized projective geometry associated
to a Jordan pair $(V^+,V^-)$,
fix $a,b \in \XX^-$ such that $\U_{a,b}:=\v_a \cap \v_b$ is non-empty,
and let 
$p: \U_{a,b} \to \XX^-$, $x \mapsto \P_{1\over 2} (a,x,b)$ be the
midpoint map assigning to $x$ the midpoint of $a$ and $b$ in
the affine space $\v_x$. Then
the set $\U_{a,b}$ becomes a symmetric space when equipped with the
multiplication map
$$
\mu(x,y)=(-1)_{x,p(x)}(y) =
\P_{-1}\bigl(x,\P_{1 \over 2}(a,x,b),y \bigr).
$$
In case $(V^+,V^-)=(J,J)$ comes from a unital Jordan algebra $J$,
then $p$ extends to a bijection denoted by
$p:\XX^+ \to \XX^-$; in this case the symmetric space $\U_{a,b}$ is
the symmetric space associated to a polar geometry (see Section 4.1), and it
is isomorphic to the quadratic prehomogeneous symmetric space $J^\times$.

\Proof. 
Let us show first that the set $\U_{a,b}$ is indeed stable under
the symmetries $\sigma_x = (-1)_{x,p(x)}$ for $x \in \U_{a,b}$.
To this end, note that, since $(-1)=(-1)^{-1}$, 
the identity (PG1) implies that, for
$x \in \U_{c,b}$, the pair 
$$
(g^+,g^-):=(\sigma_x^+,\sigma_x^-) = 
\big( (-1)_{x,p(x)}, \, (-1)_{p(x),x} \big)
$$
acts as an automorphism on $(\XX^+,\XX^-)$. 
By definition, $\sigma_x = \sigma_x^+$ is the first component of this
automorphism.
Then, using that
$\P_r(u,c,v)=\P_{1-r}(v,c,u)$ and that $r_{v,c}^{-1}=(r^{-1})_{v,c}$, we
see that
$$
\sigma_x^-(c)= \P_{-1}\big(\P_{1 \over 2}(a,x,b) , x ,a) =
\P_2(a,x,\P_{1 \over 2}(a,x,b)) = b,
$$
and similarly $\sigma_x^-(b)=a$. 
Therefore $\sigma_x^+(\v_b)=\v_{\sigma^-_x(b)}=\v_a$ and 
$\sigma_x(\v_a)=\v_b$ and hence 
$\sigma_x$ preserves $\U_{a,b}$.
Moreover, it follows that the group
$G(\U_{a,b})$ generated by all products of two symmetries preserves both
$a$ and $b$ and hence acts affinely both on 
$\v_a$ and on $\v_b$. 

The defining identies (S1) -- (S4) of a symmetric space are now all
easily checked: $x$ is fixed by $\sigma_x$ since $(-1)_{x,c}(x)=x$
for all $c \top x$; $\sigma_x^2 = \id$ since $(-1)_{x,c}^2 = 1_{x,c} = \id$;
the tangent map of $\sigma_x$ at $x$ is $-\id_{T_x M}$ since in the
chart $\v_{p(x)}$, centered at $x$, this map is just the negative of the
identity; finally, in order to prove (S3) note first that,
because of the fundamental identity (PG2), the pair $(p,p)$ is self-adjoint, 
and thus for all scalars $r$,
$$
p \circ r_{x,p(x)} (y) = p \, \P_r (x,p(x),y) = \P_r(p(x),x,p(y)) =
r_{p(x),x} \circ p(y)
$$
In particular, for $r=-1$, we have
$p \circ (-1)_{x,p(x)} = (-1)_{p(x),x} \circ p$, i.e.,
$p \circ \sigma_x = \sigma_x' \circ p$.
Using this, we obtain
$$
\eqalign{
\sigma_x \mu(y,z) & = \sigma_x \P_{-1} (y,p(y),z) = 
\P_{-1}(\sigma_x y,\sigma_x^-p(y),\sigma_x z) \cr
& = \P_{-1} (\sigma_x y,p(\sigma_x y),\sigma_x z) = 
\mu(\sigma_x y,\sigma_x z), \cr}
$$
which is precisely the property (S3).

Finally, assume that $(V^+,V^-)$ comes from a unital Jordan algebra.
In this case, there is a canonical identification of $\XX^+$ and 
$\XX^-$, denoted by $n:\XX^+ \to \XX^- $ 
(see [BeNe04], [Be03]), and then the property that $\v_a \cap
\v_b$ be non-empty is equivalent to saying that the pair
$(a,n(b))$ be quasi-invertible. If this is the case, the map $p$
is, via $n$, identified with $(-1)_{a,n(b)}$ ([Be03, 2.1]), and hence
extends canonically to a bijection $\XX^- \to \XX^+$ (called an
{\it inner polarity}). 
In [Be03, Section 6] it is shown that the symmetric space associated to
an inner polarity is, in a suitable chart, the same as the quadratic
prehomogeneous space $J^\times$.
\qed 

\msk \nin {\bf 3.2. The quasi-inverse revisited.}
In order to show that the symmetric spaces of Theorem 3.1 and of
Theorem 2.3
coincide, we need the global interpretation of the quasi-inverse
via translation groups: in any module $\m$ over a ring $\K$ 
define for $r \in \K$ the proper dilations
$r_x(y)=(1-r)x+ry$; then if $2$ is
invertible in $\K$, vector addition can be recovered from the dilations
via $x+y = 2 {x+y \over 2} = 2_0 2_x^{-1} (y)$,
i.e., translations are given by $\tau_x = 2_0 2_x^{-1}$.
In a generalized projective geometry, 
and with respect to a fixed base point $(o,o')$,
the identity (PG1) then implies that the pair
$$
(\tau_x,\tilde \tau_x)=(2_{o,o'} 2^{-1}_{x,o'}, \, \, 2^{-1}_{o',o} 2_{o',x})
$$
is an automorphism. Dually, $(\tilde \tau_a,\tau_a)$ for $a \in V^-$ is then
defined.
One can show that, in the affine chart $(V^+,V^-)$ corresponding to the
base point $(o,o')$,
$\tilde \tau_a$ is indeed the quasi-inverse from Jordan theory. Similarly, 
we have
$$
(\tau_{-x},\tilde \tau_{-x})=(2_{x,o'} 2^{-1}_{o,o'} , \, \, 
2^{-1}_{o',x} 2_{o',o}) .
$$

\Theorem 3.3. Let notation be as in {\rm Theorem 3.1} and
fix $o':=b$ as base point
in $\XX^-$ and choose an arbitrary base point $o \in \v_a \cap \v_b$.
We write $(V^+,V^-)$ for the Jordan pair associated to the base point $(o,o')$.
Then $\U_{a,b}=\U_{a,o'}$ is naturally identified with $\U_a$,  
and in the global chart $V^+$ of $\U_{a,b}$,
the symmetric space structure from {\rm Theorem 3.1} is given by the explicit
formula from {\rm Theorem 2.3}.

\Proof. We calculate in the Jordan pair $(V^+,V^-)=(\v_{o'},\v_o')$. 
For $a \in V^-$ and $x \in V^+$, the following are equivalent:
$$
\eqalign{
x \top a & \quad \iff  \quad x \in 
\v_a = \tilde \tau_a (\v_{o'})=\tilde \tau_a(V^+) \cr
& \quad \iff \quad
\tilde \tau_{-a}(x) \in V^+ \cr
&\quad \iff \quad
(x,-a) \, \hbox{quasi-invertible}.\cr}
$$
Hence
$\v_a \cap \v_{o'} = \tilde \tau_a(V^+) \cap  V^+ = \U_a$.

\ssk
Now let us show that the symmetric space structure from Theorem 3.1
is given by the explicit formulae
$$
\eqalign{
\sigma_o & = - \tilde \tau_{-a} \cr
\sigma_x \sigma_o & = \QQ(x) = \tau_{2x+Q(x)a} \circ B(x,-a). \cr}
\eqno (3.1)
$$
The proof is in several steps.
Let $d:= {1 \over 2} a$, i.e.,  $a = 2d = 2_{o,o'} d$.
Then 
$$
\eqalign{ \sigma_x(y)=
\P_{-1}(x,\P_{1 \over 2}(o',x,2d),y) & =
\P_{-1}(x,2^{-1}_{o',x} 2_{o',o}d,y)  \cr
& =\P_{-1} (x,\tilde \tau_{-x} d,y) \cr
& = \tau_{-x} \P_{-1}(2x,d,y+x) \cr
& = \tau_{-x} \circ (-1)_{2x,d} \circ \tau_x (y), \cr}
\eqno (3.2)
$$
and for $x=o$ this gives
$$
\sigma_o = (-1)_{o,d} = (-1)_{\tilde \tau_d.o,\tau_d.o'} =
\tilde \tau_d (-1)_{o,o'} \tilde \tau_{-d} = - \tilde \tau_{-2d} =
- \tilde \tau_{-a},
$$
proving the first relation from (3.1).
Now we are going to prove the second relation from (3.1).
Recall from [Be02, Cor. 5.8] or [Be00, X.3.2] that 
one can express the dilations $r_{z,b}$ for a scalar $r \in \K$
and a quasi-invertible pair $(z,b)$ via translations, quasi-inverses
and the dilation $r \id_{V^+} = r_{o,o'}$. 
For the scalar $r=-1$, the formula reads:
$$
(-1)_{z,b} =\tau_z \circ \tilde \tau_{2 \tilde \tau_{-z}(b)} \circ
\tau_z \circ (-1) = (-1) \circ
\tau_{-z} \circ \tilde \tau_{-2 \tilde \tau_{-z}(b)} \circ
\tau_{-z} 
\eqno (3.3)
$$
where $(-1)=(-1)_{o,o'}=-\id$.
Using this, we get from (3.2) along with $\sigma_o = (-1)_{o,d}$,
$$
\eqalign{
\sigma_x \sigma_o &=\tau_{-x}\circ (-1)_{2x,d} \circ \tau_x \circ (-1)_{o,d}
\cr
&  = \tau_{-x} \circ  (-1) \circ
\tau_{-2x} \circ \tilde \tau_{-2 \tilde \tau_{-2x}(d)} \circ
\tau_{-2x} \circ \tau_x \circ (-1) \circ \tilde \tau_{-2d}
\cr
& = 
 \tau_{x} \circ \tilde \tau_{2 \tilde \tau_{-2x}(d)} \circ
 \tau_x \circ \tilde \tau_{-2d}
\cr}
\eqno (3.4)
$$
From the proof of Theorem 3.1,
we know already that $\sigma_x \sigma_o$ acts as an affine transformation
on $\U_a=V^+ \cap \v_a$. Thus, in order to prove the second relation from 
(3.1), it suffices to determine its linear part and its value at $o$.
The value at $o$ of this transformation is
$$
\eqalign{
\sigma_x \sigma_o (o) &  = 
x +  \tilde \tau_{2 \tilde \tau_{-2x}(d)} (x) \cr
& = 
x + {1 \over 2} \tilde \tau_{ \tilde \tau_{-2x}(d)} (2x) \cr
& = 
x + {1 \over 2} (2x + Q(2x)d) = 
2x + 2Q(x)d = 2x + Q(x)a \cr}
$$
where we passed to the last line by using the covariance property of
the quasi-inverse map [Lo75, Th.\ I.3.7]  
which implies, in the notation of [Lo75], 
$$
\eqalign{
\tilde \tau_{\tilde \tau_{-z} b} (z) = 
z^{(b^{-z}) } & = B(-z,b) \bigl( (z-z)^b - (-z)^b \bigr) \cr
& = - B(-z,b) B(-z,b)^{-1}  (-z - Q(-z)b) = z + Q(z) b. \cr}
$$
Let us determine the linear part of $ \sigma_x \sigma_o$.
For a general automorphism $g$ such that $g(o) \in V^+$, the linear
part in the triple decomposition 
$g \in \tau_{V^+} \Aut(V^+,V^-) \tilde \tau_{V^-}$
is the inverse of the {\it denominator} $d_g(o)$ (see [BeNe04, Theorem 2.10
and its proof]). Now,
the denominator $d_g(o)$ of the transformation
$g = \tilde \tau_w \circ \tau_v$ is
$d_g(o)=B(v,w)$ (see [BeNe04]), and hence the linear part of
$\sigma_x \sigma_o$ is
$$
\eqalign{
(d_{\sigma_x \sigma_o}(o))^{-1} & =
(d_{\tilde \tau_{2 \tilde \tau_{-2x}(d)} \circ
 \tau_x} (o))^{-1}  \cr
& = B\bigl(x, 2 \tilde \tau_{-2x}(d)\bigr)^{-1} =
B\bigl(2x,  \tilde \tau_{-2x}(d)\bigr)^{-1} \cr
& = B(2x,-d) = B(x,-2d)= B(x,-a) \cr}
$$
where we used JP35 in passing to the last line.
Summing up,
$\sigma_x \sigma_0 = \tau_{2x+Q(x)a} \circ B(x,-a)$, and (3.1) is proved.
\qed

\nin Note that combining Theorem 3.1 and Theorem 3.3 gives another (in fact,
independent) proof of Theorem 2.3. Such a proof does not involve the
technique (which is artificial from a geometric point of view) of adjoining
a unit element to the Jordan algebras $V_a^+$.

\sectionheadline{4. Conformal deformations of polar geometries}

\nin  {\bf 4.1. The symmetric space of a generalized polar geometry.}
Recall from Appendix B.3 that Jordan triple systems $(V,T)$ 
correspond to {\it generalized polar geometries} $(\XX^+,\XX^-;p)$, i.e.,
to generalized projective geometries with polarity $p:\XX^+ \to \XX^-$.
We denote by
$$
M^{(p)} := \{ x \in \XX^+ | \, x \top p(x) \}
\eqno{\rm (4.1)}
$$
the set of non-isotropic points of the polarity $p$. This set becomes a
symmetric space when equipped with the product
$$
\mu(x,y) :=  \sigma_x(y):= (-1)_{x,p(x)}(y)=\P_{-1}\bigl( x,p(x),y \bigr)
\eqno{\rm (4.2)}
$$
(see [Be02, th.\ 4.1]; the proof is similar to the one
of Theorem 3.1). 
Fixing a point $o \in M^{(p)}$ and choosing the transversal pair
$(o,o')=(o,p(o))$ as base point in $(\XX^+,\XX^-)$, we can describe the
affine image of $M^{(p)}$ as follows: 
when we identify $\XX^+$ and $\XX^-$ and $V^+$ and $V^-$ 
via $p$, then $(x,p(x))$ is a transversal
pair  iff $(x,-x)$ is a quasi-invertible pair (see Chapter 3),
i.e., iff the Bergman operator $B(x,-x)$ is invertible, 
and hence\footnote{$^1$}{\eightrm
One should note that the 
sign in the correspondence between spaces and triple systems
is a matter of convention. 
Geometers look at Hermitian symmetric spaces of non-compact type 
as having ``negative"
curvature, and on the other hand, one is used to associate them to 
``positive" Hermitian JTS, see [Lo77]. 
In order to check signs, one may use that
for the JTS $\R$ with $T(x,y,z)=2xyz$, the
positive Bergmann polynomial $B(x,-x)=(1+x^2)^2$ must belong to the compact 
space $\R \PPP^1$ and the ``opposite'' Bergmann polynomial 
$B(x,x)=(1-x^2)^2$ to its non-compact dual.}
$$
M^{(p)} \cap V = \{ x \in V | \, B(x,-x) \, {\rm is \, invertible} \, \}.
\eqno{\rm (4.3)}
$$
%
%(Note that this affine image is of a different kind than the one considered
%in Chapter 2. It corresponds to the ``centered'' vectorialization with
%respect to the pair $(o^+,o^-)=(o,p(o))$, whereas in Theorem 3.1 $o'$
%was different from $p(o)$.) 
%

\msk
\nin  {\bf 4.2. The symmetric space associated to the $\alpha$-homotope.}
Now consider the Jordan
triple systems $V_\alpha = (V,T_\alpha)$, the $\alpha$-homotopes of $V$
with respect to structural  transformations $\alpha:V \to V$ (cf.\ Lemma 1.12).
Since $V_\alpha$ is a JTS in its own right, the preceding construction
associates to $V_\alpha$ a 
generalized polar geometry $(\XX_\alpha^+,\XX_\alpha^-,p_\alpha)$
with symmetric space $M_\alpha$. 
The explicit formula for the affine image (analog of the ``bounded 
realization") is of the same type as (4.3),
where the underlying $\K$-module is the same, but $T$ is replaced by
$T_\alpha$ and $B$ by $B_\alpha(u,v)=B(u,\alpha v)$.
However, {\sl a priori},
the global spaces $\XX_\alpha^+$ for various $\alpha$, are different spaces.
 In order to compare the geometries
$\XX_\alpha$ and the symmetric
spaces $M_\alpha$ with each other and to realize them as deformations
of $M^{(p)}$, we need a common realization, which
is given by the following result:

\Theorem 4.3. 
Let $(\XX^+,\XX^-;p)$ be the generalized polar geometry  of 
a Jordan triple system $(V,T)$, and let $\alpha:V \to V$ be a
structural transformation, i.e., $\alpha \in \Svar(T)$. 
Let $V_\alpha=(V,T_\alpha)$ be the $\alpha$-homotope of $(V,T)$ and 
let $(\XX_\alpha^+,\XX_\alpha^-;p_\alpha)$ 
be the generalized polar geometry associated to  $(V,T_\alpha)$, 
with symmetric space denoted by $M_\alpha \subset \XX_\alpha^+$.
Then there exists a canonical embedding
$\Phi_\alpha:\XX_\alpha^+ \to \XX^+$ having the following properties:
\ssk
\item{(1)} 
The affine image of $M_\alpha$ under 
this imbedding coincides with the description of
$M_\alpha$ in the chart $V_\alpha$, i.e.
$$
\eqalign{
\Phi_\alpha(M_\alpha) \cap V 
% & = \{ x \in V \vert \, x \top p \alpha(x) \} \cr 
& =  \{ x \in V \vert \, B(x,-\alpha x) \in \Gl(V) \}      \cr}
$$
\item{(2)} 
The tangent map of $\Phi_\alpha$  at
the origin can be algebraically defined, and it is the identity map
$V_\alpha \to V$.
\item{(3)}
If $\alpha$ is invertible, then the image of this imbedding
is precisely the symmetric space $M^{(p \alpha)}$ associated to
the polarity $p \circ \alpha : \XX^+ \to \XX^-$.

\Proof. First we describe the construction of $\Phi_\alpha$. Saying that
$\alpha:V \to V$ is a structural transformation of the JTS $(V,T)$
is equivalent to saying that the map of Jordan pairs
$$
(\alpha,\alpha):(V,V) \to (V,V)
$$
is a structural transformation of the Jordan pair $(V,V)$
 (see Appendix C). 
We apply Lemma C.3, with $(W^+,W^-)=(V^+,V^-)$ and $(f,g)=(\alpha,\alpha)$;
then the Jordan pair $(V^+,W^-)$ 
from Lemma C.3 is $(V_\alpha,V_\alpha)$, with 
the JTS $V_\alpha$ defined by Lemma 1.12.
The second point of Lemma C.3 says that
$$
\phi=(\id , \alpha): (V_\alpha , V_\alpha) \to (V , V)
$$
is a homomorphism of Jordan pairs. 
%(It is not a JTS-homomorphism, unless
%$\alpha = \id$, because it does not commute with the exchange automorphisms!)
As in Section C.4,
the homomorphism $(\id,\alpha)$ induces a homomorphism
$$
(\Phi^+,\Phi^-): (\XX^+(V_\alpha),\XX^-(V_\alpha)) \to (\XX^+(V),\XX^-(V))
$$
with injective first component 
$\Phi^+:\XX^+(V_\alpha) \to \XX^+(V)$.

We define $\Phi_\alpha:= \Phi^+$ and show
that this imbedding has the desired properties.
Since $(\Phi^+,\Phi^-)$ is a homomorphism of generalized projective
geometries, its restriction to affine parts with respect to the
base points is given by the pair of linear maps which induces it, 
namely by $(\id,\alpha)$, and hence
the algebraic tangent map can be defined and is equal to
the pair $(\id,\alpha)$ 
(cf.\ [Be02]).  Thus the tangent map of $\Phi^+$ 
at the base point $o$ is the identity map $V_\alpha \to V$.
This proves (2).

We prove (1). Let $x \in V^+=V$. 
Recall that the element $p(x) \in \XX^-$ is identified with $x \in
V^- =V$.
By definition, $x$ belongs to the symmetric space $M^{(p)}$
if and only if $x$ and $p(x)$ are transversal, if and only if
the pair $(x,-p(x))=(x,-x)$ is transversal (see proof of 3.1), if and only if
$B(x,-x)$ is invertible.
The same argument, applied to the JTS $T_\alpha$, shows that $x$ belongs
to $M_\alpha$ if and only if $B_\alpha(x,-x)$ is invertible. But 
$B_\alpha(x,y) = B(x,\alpha y)$ (Lemma 1.12),
and hence $x$ belongs to $M_\alpha$ if and only if $B(x,-\alpha x)$ is
invertible.

We prove (3).  Assume
$\alpha$ is $T$-selfadjoint and invertible.  Then 
$(\id,\alpha):(V_\alpha,V_\alpha) \to (V,V)$ is a Jordan pair isomorphism,
and the induced map $(\Phi^+,\Phi^-)$ is an isomorphism of geometries.
Under this isomorphism, the polarity of $(\XX^+(V_\alpha),\XX^-(V_\alpha))$
defining the symmetric space $M_\alpha$ (the exchange automorphism)
corresponds to the map $p \circ \alpha$ which is therefore the polarity
defining the symmetric space $\Phi_\alpha(M_\alpha)$.
This had to be shown.
\qed

\Remark 4.4. (Realisation as a graph.)
Assume that $\alpha$ is a $T$-symmetric map. Then
it is easily checked that the graph of $\alpha$,
$\Gamma_\alpha = \{ (v,\alpha v) \vert \, v  \in V \}$,
is a sub-JTS of the polarized JTS $(V \oplus V,\tilde T)$ with
$$
\tilde T\big( (x,a),(y,b),(z,c) \big) = (T(x,b,z),T(a,y,c)),
$$
and that $\Gamma_\alpha$ is isomorphic to $V_\alpha$ as JTS.
By the preceding construction, this situation globalizes as follows:
the polarized JTS $(V \oplus V,\tilde T)$ corresponds to the 
polar geometry $(\XX^+ \times \XX^-,\XX^- \times \XX^+)$ with the
exchange automorphism (see [Be02, Prop.\ 3.6]) and with associated
symmetric space
$M=\{ (x,\alpha) \in \XX^+ \times \XX^- \vert \, x \top \alpha \}$. Then
$M_\alpha$ is the intersection of $M$ with the graph of the map 
$$
 \Phi^+(\XX^+) \to \XX^-, \quad \Phi^+(x) \mapsto \Phi^-(x).
$$

\ssk
\Remark 4.5. (Alternative proof in the real finite-dimensional case.) 
In the finite-dimensional real or complex case, the main
statements of Theorem 4.3 can be proved in a different way:
if $\g = \g_1 \oplus \g_0 \oplus \g_{-1}$ with involution $\theta$ 
is the involutive $3$-graded Lie algebra associated to the JTS $T$, then
$$
\q_\alpha := \{ v + \theta(\alpha v) \vert \, v \in \g_1 \}
$$
is a sub-Lie triple system of $\g$ (a ``defining plane'' in the terminology
of [Ri70]). Let $\h_\alpha = [\q_\alpha,\q_\alpha]$ and
define the subalgebra $\g_\alpha := \h_\alpha \oplus \q_\alpha$ of $\g$.
Let $G_\alpha$ and $H_\alpha$ be the analytic subgroups of $G = \Aut(\g)$
with Lie algebra $\g_\alpha$, resp.\ $\h_\alpha$.
Then the orbit $M_\alpha := G_\alpha .o$ of the base point $o$ in $\XX^+$
is open, and it is a homogeneous symmetric space isomorphic to $G_\alpha/
H_\alpha$ coinciding (possibly up to topological connected components) 
with the space $M_\alpha$ from Theorem 4.3.
This construction can already been found in the work [Ri70] by
A.A.\ Rivillis. Since it uses the exponential map as essential ingredient,
it gives rather results of local nature, and it does not work for general
base fields.

\Theorem 4.6.
Assume that $\K$ is a topological ring with dense unit group, and that
$(V,T)$ 
is a topological Jordan triple symstem over $\K$ such that
\item{(C1)}
The set $(V \times V)^\times = \{(x,y) \in V \times V \vert \, B(x,y) \in \Gl(V) \}$ of 
quasi-invertible pairs is open in $V \times V$, and the map
$$
(V \times V)^\times \times V \to V, \quad 
(x,y,z) \mapsto B(x,y)^{-1}z
$$
is continuous.

\ssk
\nin Assume moreover that $\alpha$ is a continuous structural transformation
 map.
Then $M_\alpha$ is a smooth symmetric space over $\K$, and it is an open
submanifold of $\XX=\XX^+$. The Lie triple system of the pointed symmetric 
space $(M_\alpha,o)$ is $V$ equipped with the Lie triple product
$$
[X,Y,Z]_\alpha :=  T_\alpha(X,Y,Z) - T_\alpha(Y,X,Z) =
T(X,\alpha Y,Z) - T(Y,\alpha X,Z).
$$

\Proof.
This follows directly from [BeNe05, Theorem 6.3] since
the JTS $T_\alpha$ again satisfies the property (C1).
\qed

\nin
In particular, the LTS $[\cdot,\cdot,\cdot]_\alpha$ and
$[\cdot,\cdot,\cdot]_{-\alpha}$ are negatives of each other, which means
that the symmetric spaces $M_\alpha$ and $M_{-\alpha}$ are
{\it c-duals} of each other. 

\sectionheadline{5. Behavior of geometric structures under deformation}

\nin
For applications in geometry and harmonic analysis one not only needs to
know how the  affine image $M_\alpha \cap V$
of the symmetric space $M_\alpha$ as a set, but also how its
geometric structure  depends on the deformation parameter $\alpha$.
We give some answers to this question. As a rule, the dependence on the
 parameter $\alpha$
can be made explicit if we simply mind the rules: 
in the following formulae, replace 
$M$ by $M_\alpha$,
$T$ by $T_\alpha$,
$B(x,y)$ by $B(x,\alpha y)$ and $Q(x)$ by $Q(x) \circ \alpha$.

\msk \nin
{\bf 5.1. The product map.}
The product map of the symmetric space $M$ is given by Formula (3.3), 
by letting $x=z$, $b=p(x)$:
$$
\eqalign{
\sigma_x(y) = (-1)_{x,p(x)}(y)& 
% = 
% \tau_x \tilde \tau_{\tilde \tau_{-x} (x)} \circ (-1)_{o,o'} \circ
% \tilde \tau_{-\tilde \tau_{-x} (x)} \tau_{-x} (y)
% \cr
%  &=(-1)_{o,o'} \tau_{-x} \tilde \tau_{-2\tilde \tau_{-x} (x)}  \tau_{-x} (y)
% \cr 
= 
 \tau_{x} \tilde \tau_{2\tilde \tau_{-x} (x)}  \tau_{x} (-1)_{o,o'}(y) \cr 
&= 
\tau_{2 \tilde \tau_x (x)} B(2 \tilde \tau_x(x),-x)^{-1} 
\tilde \tau_{2 \tilde \tau_x (x)} (-y)
\cr}
\eqno (5.1)
$$
In particular, the ``square'' $x^2 := \sigma_x(o)$ in the symmetric space $M$
is
$$
\sigma_x(o)=2 \tilde \tau_x(x) = 2 (\id_V - Q(x))^{-1} x
\eqno (5.2)
$$
(the last equality is proved as in [Be00, p. 195]).  

\msk \nin {\bf 5.2. Second derivatives and second order tangent bundles.}
The tangent bundle $TM$ of $M$, which is again a symmetric space, is simply
obtained by scalar extension by dual numbers, as already used in Section 2.1
(cf.\ [BeNe05, Th.\ 5.5]). Hence the geometrically important second tangent
bundle $TTM$ is given by  
scalar extension by the second order tangent ring
$TT\K = \K[\eps_1,\eps_2]$. The same statement holds for the geometry
$(\XX^+,\XX^-)$ itself. For instance, by a straightforward calculation 
using that $\eps_1^2=0=\eps_2^2$, we get from the definition of the
quasi-inverse and from (5.1)
$$
\eqalign{
\tilde \tau_a(\eps_1 v_1 + \eps_2 v_2) & = B(\eps_1 v_1 + \eps_2 v_2,a)^{-1}
(\eps_1 v_1 + \eps_2 v_2 - Q(\eps_1 v_1 + \eps_2 v_2)a) \cr
& = \eps_1 v_1 + \eps_2 v_2 + \eps_1 \eps_2 T(v_1,a,v_2), \cr
\sigma_x(x + \eps_1 v_1 + \eps_2 v_2) & = x - \eps_1 v_1 -
\eps_2 v_2 + \eps_1 \eps_2   T(v_1,\tilde \tau_{-x}(x),v_2). \cr}
$$
On the other hand, ordinary second differentials and second order tangent maps
of twice continuously differentiable maps $f$ are related via 
$$
f(x+ \eps_1 v_1 + \eps_2 v_2) = f(x) +
\eps_1 df(x)v_1 + \eps_2 df(x) v_2 + \eps_1 \eps_2 d^2 f(x)(v_1,v_2)
$$
(see [Be06, Eqn.\ (7.18)]), so that from the preceding calculation we get
the second differential
of $\tilde \tau_a$ at the origin and the one of $\sigma_x$ at $x$.

\msk \nin {\bf 5.3. The canonical affine connection.}
It is possible to define the canonical affine connection of the symmetric
space in a purely algebraic way. In the situation of Theorem 4.6 it 
coincides with any of the various differential geometric definitions of 
a connection (cf.\ [Be06, Chapter 26]). 
More precisely,
in [Be06, Chapter 26] it is shown that the canonical connection $\nabla$
 of a symmetric
space $M$ can be defined in terms of the second order tangent bundle $TTM$,
and we have the following chart description: let $x \in M$, $V$ a chart
around $x$ inducing the flat connection $\nabla^0$ on some neighborhood
of $x$; then the difference $C:=\nabla - \nabla^0$ (i.e., the
Christoffel tensor of the canonical connection with respect to the chart)
at the point $x$ is simply the second
derivative of the symmetry $\sigma_x$ at $x$:
$C_x = d^2 \sigma_x (x)$. Since we have calculated this second differential
above, we get
$$
C_x (u,v) =  T(v_1,\tilde \tau_{-x}(x),v_2)=T(v_1,(\id_V - Q(x))^{-1} x ,v_2)
\eqno (5.3)
$$
(for the last equality, see (5.2)).
In the finite-dimensional real case this result is already given in
[BeHi01, Theorem 2.3]; however, the proof given there does not carry over
to the general situation considered here.
As explained in [BeHi01], (5.3) implies
that the canonical connections of the symmetric spaces $M_\alpha$ are
{\it conformally equivalent} among each other.
Note also that, for $\alpha=0$, we get $C_x=0$ for all $x$, 
and hence $\nabla=\nabla_0$, as expected.

\msk \nin {\bf 5.4.  Trace form, pseudo-metrics and their duals.}
It is clear that semi-simplicity is not preserved under deformations.
Let us assume that
$\K$ is a field, the JTS $V$ is finite-dimensional over $\K$, and that
 the {\it trace form\/}
 $g_0(u,v) := \tr (T(u,v,\cdot))$ is non-degenerate.
If this is the case, then $g_0$ extends to a pseudo-Riemannian metric tensor
(non-degenerate symmetric
two fold covariant tensor field) with affine picture given by
$$
g_x(u,v) :=g_o(u,B(x,-x)^{-1} v), \quad (u,v \in V, \, x \in M)
\eqno (5.4)
$$
(see [Be00, X.6] for the case $\K=\R$).
On the other hand, the dual two fold contraviant tensor field with affine
picture $\gamma_x(\phi,\psi)=g_o(\phi,B(x,-x)\psi)$ can always be defined
and exists also for the deformed spaces.

\msk \nin {\bf 5.5. Invariant measures: group case.}
Let $\K=\R$ and $A$ a finite dimensional associative unital $\R$-algebra. 
In general, the Lie group $G=A^\times$ and its deformations $G_a$
will not be unimodular (take, e.\ g., the algebra of upper triangular
square matrices). However, if $G=A^\times$ happens to be unimodular
(e.g, the case of the general linear groups), then also its deformations
$G_a$ remain unimodular:
the modular function is $g \mapsto \det (\Ad(g))$;
if this function is trivial for all 
groups $G_a$, $a \in A^\times$, then by density of $A^\times$ in $A$
it remains trivial for all $a \in A$.

\msk \nin {\bf 5.6. Invariant measures on symmetric spaces.}
Let $\K=\R$ and $V$ a finite-dimensional JTS over $\R$.
Then in general there is no measure on $M$ which is invariant under the
transvection group $G(M)$. Such a measure exists if $(V,T_\alpha)$ 
has non-degenerate trace form:
in this case, the corresponding invariant integral is given by the 
explicit formula
$$
I_\alpha(f) = \int_{M_\alpha \cap V} 
f(x) \vert \det B(x,-\alpha x) \vert^{-{1\over 2}} \,dx ,
\eqno (5.5)
$$
where $dx$ is  Lebesgue measure on $V$ (see [Be00, X.6.2]).
The term $ \vert \det B(x,-\alpha x) \vert^{-{1\over 2}} \,dx $
translates the existence of an invariant half-density on the symmetric
space $M_\alpha$ (cf.\ [BeHi00]). But one may define a half-density
on $M_\alpha \cap V$ by the same formula for any $\alpha$ in the
structure variety.
Since the invariance condition is an equality of two functions
depending continuously on $\alpha$, these equality will still hold
for all $\alpha$ in the closure of the subset $S$ of the structure
variety for which the trace form of $T_\alpha$ is non-degenerate,
and hence for $\alpha \in \overline S$ there still exists an 
invariant measure on $M_\alpha$.
(It is not clear wether this property holds for all $\alpha$ in the
structure variety. There may be topological
 connected compenents of the structure
variety containing only ``degenerate'' $\alpha$'s, see end of Section 6.5.)
 In particular, for $\alpha \to 0$, the invariant measure tends to the Lebesgue
measure of $V$. 

\msk \nin {\bf 5.7. Further results and topics.}
Further results on geometric structures on the spaces $M_\alpha$ in the
real, non-degenerate case can be found in the paper
[Ma79] by B.O.\ Makarevic, and it would certainly be interesting to generalize
some of these results.
Finally, one would like to study the dependence of 
harmonic analysis and representation theory
on the symmetric spaces $M_\alpha$ 
on the deformation parameter $\alpha$.
Although the situation here is different from the one considered, e.g.,
in [DR85], it seems possible to transfer ideas developed in the context
of contractions of semisimple Lie groups to the set-up just described.
Since all relevant geometric data of $M_\alpha$ 
depend algebraically or at least analytically on $\alpha$, one may
expect that something similar is true on the level of data of harmonic
analysis.
See, for instance [FaPe04] and [Pe02] for examples of  harmonic analysis
on symmetric spaces of the type $M_\alpha$.
The parameter $\alpha$ is fixed in these works, 
so that the dependence on $\alpha$ is not
visible. When varying $\alpha$, not only the aspect of deformations of $\alpha$
to singular values is of interest, but also the one of {\it analytic
continuation} between isotopic, but non isomorphic values of $\alpha$.
For instance, the duality between {\it compactly causal} and 
{\it non-compactly causal symmetric spaces} is reflected by the isotopy
between $\alpha$ and $-\alpha$ (cf.\ [Be00, XI.3]).

\sectionheadline
{6. Some examples}

\msk \nin {\bf 6.1. Quadratic matrices and
associative algebras.} Let us return to the setting of Section 0.2 and
consider the associative algebra
$A = M(n,n;\K)$  of square matrices over a field $\K$. 
Then the $A^\times \times A^\times$-orbits are represented by the
idempotent matrices
$$
e:= e_r := \pmatrix{\1_r & 0 \cr 0 & 0 \cr} , \quad \quad r=0,\ldots,n.
$$
The group $G_e=\{ x \in A \vert \, \1 + ex \in A^\times \}$ 
is explicitly given by
$$
G_e = \Big\{ x=\pmatrix{\alpha & \beta \cr \gamma & \delta \cr} \vert \,
\det \pmatrix{\1_r + \alpha & \beta \cr 0 & \1_{n-r} \cr} \not= 0 \Big\} =
 \Big\{ \pmatrix{\1_r + \alpha & \beta \cr \gamma & \delta \cr} \vert \,
\det(\1_r + \alpha) \not= 0 \Big\}
$$
(where $\beta,\gamma,\delta$ are arbitrary),
with product $x \diamond_e y  = xey+x+y$, i.e.
$$
\pmatrix{\alpha & \beta \cr \gamma & \delta \cr}
\diamond_e
\pmatrix{\alpha' & \beta' \cr \gamma' & \delta' \cr}  =
\pmatrix{\alpha \alpha' & \alpha \beta' \cr
\gamma \alpha' & \gamma \beta'\cr}
+\pmatrix{\alpha & \beta \cr \gamma & \delta \cr}
+\pmatrix{\alpha' & \beta' \cr \gamma' & \delta' \cr}  .
$$
This group is isomorphic to a semidirect product of $\Gl(r,\K)$ and
a Heisenberg type group: 
for $\beta,\gamma,\delta$ zero, we get a subgroup isomorphic to $\Gl(r,\K)$,
and for $\alpha=0$, we obtain a subgroup with product
$$
\pmatrix{0 & \beta \cr \gamma & \delta \cr}
\diamond_e
\pmatrix{0 & \beta' \cr \gamma' & \delta' \cr}  =
\pmatrix{0 & 0 \cr 0 & \gamma \beta'\cr}
+\pmatrix{0 & \beta \cr \gamma & \delta \cr}
+\pmatrix{0 & \beta' \cr \gamma' & \beta' \cr}  .
$$
It is easily checked that $G_e$ is the semidirect product of these
two subgroups. The Lie algebra of $G_e$ is given by the bracket
$[x,y]_e = xey-yex$, which yields explicitly
$$
\Big[ \pmatrix{\alpha & \beta \cr \gamma & \delta \cr} ,
\pmatrix{\alpha' & \beta' \cr \gamma' & \delta' \cr} \Big]_e =
\pmatrix{\alpha \alpha'- \alpha'\alpha & \alpha \beta'- \alpha' \beta 
\cr\gamma \alpha' - \gamma' \alpha & \gamma \beta' - \gamma' \beta \cr}
$$
so that for $\alpha = 0 = \alpha'$ we get the 2-step nilpotent algebra
$$
\Big[ \pmatrix{ 0 & \beta \cr \gamma & \delta \cr} ,
\pmatrix{ 0 & \beta' \cr \gamma' & \delta' \cr} \Big]_e =
\pmatrix{ 0 & 0 \cr 0 & \gamma \beta' - \gamma' \beta \cr}.
$$
If $A$ is a general associative algebra and $e$ an idempotent in $A$, 
then the preceding calculations
carry over, where matrices have to be understood with respect to
the {\it Peirce-decomposition} of $A$, where $c:=\1-e$: 
$$
A = cA \oplus eA = cAc \oplus cAe \oplus eAc \oplus eAe.
$$
Then $G_e$ is isomorphic to a semidirect product of 
$(eAe)^\times$ with a Heisenberg type group $cAe \oplus eAc \oplus cAc$.
However, in general not every element of $A$ is similar to an idempotent,
and hence other types of deformations can appear.

\msk
\nin
{\bf 6.2. Rectangular matrices.}
We consider the Jordan pair $(V^+,V^-)=(M(p,q;\K),M(q,p;\K))$ 
of rectangular matrices over a field $\K$, with trilinear maps
$T^\pm(x,y,z)=xyz+zyx$.  The group
$\Gl(p,\K) \times \Gl(q,\K)$ acts by automorphisms of this Jordan
pair(in the usual way from
the right and from the left), and every matrix $a \in V^-$ is conjugate
under this action to a matrix of the form
$e=e_r$ with $e_{ii}=1$ for $i=1,\ldots,r$ and $0$ else.
Thus we may assume $a=e$, and then the symmetric space structure on
$\U_e$ from Theorem 2.3 is calculated in a similar way as above.
Recall that a pair $(x,y) \in (M(p,q;\K),M(q,p;\K))$ is quasi-invertible
if and only if the matrix $\1_p - xy$ is invertible in $M(p,p;\K)$;
therefore $\U_e$ is the set of matrices $x \in M(p,q;\K)$ such that
$\1_p + xe$ is invertible, i.e.
$$
\U_e = 
\Big\{ \pmatrix{\1_r + \alpha & \beta \cr \gamma & \delta \cr} \vert \,
\det(\1_r + \alpha) \not= 0 , \quad
\beta,\gamma,\delta \, {\rm arbitrary} \,\Big\}.
$$
Using that $Q(x)e=xex$ and $B(x,-e)y=(\1+xe)y(\1+ex)$, we get the
quadratic map of the symmetric space $\U_e$
$$
\QQ(x)y = 2x + Q(x)e + B(x,-e)y=2x + y +xex  + xey + yex + xeyex.
$$
Again one finds that $\U_e$ has a bundle structure: the basis is given
be letting $\beta,\gamma,\delta$ equal to zero; it is the space
of $r \times r$-Matrices, with the symmetric space structure
$$
\QQ(\alpha)\alpha'=2 \alpha + \alpha' + \alpha^2 +
\alpha \alpha' + \alpha' \alpha +
\alpha \alpha' \alpha 
$$
which is isomorphic to the group $\Gl(r,\K)$, seen as a symmetric space
(group case, cf.\ Appendix A.3). The fiber is given by letting $\alpha=0$, 
and we get explicitly 
$$
\QQ\Big( \pmatrix{0 & \beta \cr \gamma & \delta \cr} \Big)
\pmatrix{0 & \beta' \cr \gamma' & \delta' \cr}  =
\pmatrix{0 & 2 \beta + \beta' \cr
2 \gamma +\gamma' & 2 \delta + \delta' \cr} + 
\pmatrix{0 & 0 \cr 0 &  \gamma \beta + \gamma \beta' + \gamma' \beta \cr} .
$$
Knowing that the quasi-inverse in the Jordan pair of rectangular matrices
is given by $\tilde \tau_a (x) = (\1-xa)^{-1} x = x (\1- ax)^{-1}$,
we find that the symmetry with respect to the origin in the fiber is
$$
\sigma_0 \Big(
\pmatrix{0 & \beta \cr \gamma & \delta \cr} \Big) =
\pmatrix{0 & -\beta \cr -\gamma & -\delta + \gamma \beta \cr}, 
$$
so that the product map of the fiber is finally
$$
\mu \Big( \pmatrix{0 & \beta \cr \gamma & \delta \cr},
\pmatrix{0 & \beta' \cr \gamma' & \delta' \cr} \Big) =
\pmatrix{0 & 2 \beta - \beta' \cr
2 \gamma -\gamma' & 2 \delta - \delta' \cr} + 
\pmatrix{0 & 0 \cr 0 &  (\gamma -\gamma')(\beta- \beta') \cr} .
$$
In contrast to the preceding section, where the fiber was a two-step
nilpotent group, for the symmetric space structure the fiber is now
simply a flat a symmetric space. This is not surprising since the
Lie triple  bracket on a two-step nilpotent group is zero. 
Explicitly, one can check by direct calculation that the map
$$
F : \pmatrix{0 & \beta \cr \gamma & \delta \cr} \mapsto
\pmatrix{0 & \beta \cr \gamma & \delta + \gamma \beta \cr} 
$$
is an isomorphism from the fiber onto the same set, but equipped with
the flat symmetric space structure $2x-y$. 
(The situation is completely similar to the one considered in [Be06],
where the symmetric space structure on the second order tangent bundle
$TTM$ of a symmetric space $M$ is calculated, see loc.\ cit., Chapters 26
and BA.)

\msk \nin
{\bf 6.3. Symmetric spaces associated to idempotents.}
The preceding calculations and remarks can be generalized for an
arbitrary Jordan pair $(V^+,V^-)$, as follows:
we assume that $e:=e^- \in V^-$ is an element which can be completed to an
idempotent, i.e., there existes $e^+ \in V^-$ such that
$Q^+(e^+)e^-=e^+$ and $Q^-(e^-)e^+=e^+$. (In other words, $e$ is
{\it von Neumann regular}, cf.\ [Lo75].) Associated to an idempotent
$(e^+,e^-)$, there is a
{\it Peirce-decomposition} $V^+=V_2^+ \oplus V_1^+ \oplus V_0^+$ which is the
eigenspace decomposition of the operator $T(e^+,e^-)$.
The space $V_1^+ \oplus V_0^+$ depends only on $e^-$ since
$V_1^+ \oplus V_0^+ = \ker Q(e^-)$, and $V_2^+$ depends only on $e^+$ since
$V_2 = Q(e^+)V^-$. Moreover, $V_2^+$ is a unital Jordan algebra with unit
element $e^+$ (see [Lo75] for all this).

\Theorem 6.4. {\rm (Bundle structure of $\U_e$.)}
Assume $e = e^- \in V^-$ as a von Neumann regular element and
choose $e^+ \in V^+$ such that $(e^+,e^-)$ is an idempotent.
\item{(1)}
$\U_e = \{ x = x_0+x_1+x_2 \vert \, x_0 \in V_0^+,
x_1 \in V_1^+,  x_2 + e^+ \in (V_2^+)^\times  \,  \}$
\item{(2)}
The space $\U_e \cap V_2^+ = e^+ - V_2^\times$ is a symmetric subspace
of $\U_e$ (called the {\rm base}), isomorphic to the quadratic
pre-homogeneous symmetric space $(V_2^+)^\times$.
\item{(3)}
The space $V_1^+ \oplus V_0^+$ is a symmetric subspace of $\U_e$ (called the
{\rm fiber}), and its symmetric space structure is 
isomorphic to the flat symmetric structure $2x-y$ on this set.

\Proof.
(1) We use the identity JP24 of [Lo75] which gives
$$
B(x,e)B(x,-e)=B(x,Q(e)x)=B(Q(x)e,e).
$$
Assume $x \in V_0^+ \oplus V_1^+$, i.e., $Q(e)x=0$. Hence
$B(x,e)B(x,-e)=B(x,0)=\id$. In the same way, $B(x,-e)B(x,e)=\id$, and
so $B(x,-e)$ is always invertible, with inverse $B(x,e)$. This shows that
$V_0^+ \oplus V_1^+ \subset \U_e$.
Moreover, then $\tilde \tau_x (e)=e$ since
$B(e,x)^{-1}(e-Q(e)x)=B(e,-x)(e-Q(e)x)=B(e,-x)e=e$.

Assume $x \in V_2^+$. Then $B(x,-e)$ is invertible in $\End(V^+)$
 if and only if
its restriction to $V_2^+$ is invertible, and hence, by the results on
the Jordan algebra case (Prop.\ 2.2), we see that this is equivalent to saying 
that $e^+ + x$ is invertible in $V_2^+$.

Now let $z \in V^+$ and decompose
 $z=x+y$ with $y \in V_0^+ \oplus V_1^+$ and 
$y \in V_2^+$. Then $(y,e)$ is quasi-invertible with
$\tilde \tau_y(e)=e$, as we have seen above, and
hence we may apply the identity JP34 of [Lo75] in order to obtain 
$B(x+y,e)= B(x,\tilde \tau_y (e))B(y,e) = B(x,e)B(y,e)$.
This operator is invertible in $\End(V^+)$
if and only if $x$ and $y$ satisfy the previously described
conditions, and this proves the claim.

(2)
This follows directly from (1) since the formula for the product
$\mu(x,y)$ with $x,y \in (V_2^\times +e^+)$ uses only data which
depend on the Jordan algebra $V_2^+$.

(3) We first calculate  $\sigma_0(y)$ for an element
$y \in \ker Q(e)$: using $e^{-y}=e$, we get from the
symmetry formula ([Lo75, 3.3])
$$
\eqalign{\sigma_0(y) & =
\tilde \tau_{e}(-y)  = (-y)^e = - y + Q(y) \, e^{-y} =
= - y + Q(y)e .
\cr}
$$
Next we calculate the map $\QQ(x)y=2x+Q(x)e +  B(x,-e)(y)$.
Since
$B(x,-e)y = y + T(x,e,y)+Q(x)Q(e)y=y+T(x,e,y)$,
we get
$$
\QQ(x)y=2x+y+Q(x)e+T(x,e,y) = 2x+y + {1 \over 2} T(x,e,x) + T(x,e,y).
$$
As in the matrix case (Section 6.2), it is checked by a direct calculation
that the map
$$
F: V_0^+ \oplus V_1^+ \to V_0^+ \oplus V_1^+, \quad z \mapsto z + Q(z)e
$$
is an isomorphism from the symmetric space $V_0^+ \oplus V_1^+$ with the
usual flat structure $(x,y) \mapsto 2x-y$ 
onto the symmetric space $V_0^+ \oplus V_1^+$ equipped with the symmetric 
space structure that we have just described.
\qed

\nin
The preceding proof was purely Jordan-theoretic.
It would be interesting to state and prove Theorem 6.4 in a geometric way,
i.e., in the context of Theorem 3.1. The starting point should be a
suitable geometric interpretation of von Neumann regularity. This is
indeed possible and will be taken up elsewhere.

\msk \nin {\bf 6.5. Grassmannians.}
Next we give some examples illustrating Theorem 4.3.
Assume $W=\K^{p+q}$ is a vector space over a field $\K$ of characteristic 
different from $2$ and 
let $(\XX^+,\XX^-)=(\Gras_p(W),\Gras_q(W))$ the Grassmannian geometry 
of $p$- and  $q$-dimensional subspaces of $W$, respectively.
 Then $(\XX^+,\XX^-)$ is a generalized
projective geometry; with respect to the base point $(o^+,o^-)= 
(\K^p \oplus 0,0 \oplus \K^q)$, it corresponds 
to the Jordan pair $(M(q,p;\K),M(p,q;\K))$.
We define a family of symmetric space structures associated to 
$(\XX^+,\XX^-)$:
assume $\beta:W \times W \to \K$ is a symmetric or 
skew-symmetric bilinear form on $W$ such that $\beta$ is non-degenerate
on the first summand $o^+=\K^p$. If $\beta$ is non-degenerate on all of
$W$, then $p:\XX^+ \to \XX^-$,
$E \mapsto E^\beta$,
where $E^\beta$ is the orthogonal complement of $E$ with respect to $\beta$,
is a globally defined polarity, but if $\beta$ is degenerate, 
then $p$ is only defined on the 
% (open dense)
set $\{ E \in \XX^+ \vert \, \dim (E^\beta) = q \}$.
% (in other words, 
% $W \to E^*$, $w \mapsto \beta(w,\cdot)\vert_E$ has to be surjective).
In any case,
the space of non-degenerate linear subspaces of $W$, 
$$
M^\beta:= \{ E \subset W | \, W = E \oplus E^\beta \}
$$
is a symmetric space with product map
$$
\mu(E,F)=\sigma_E (F) = (-1)_{E,E^\beta} (F),
$$
where, for any scalar $r \in \K$ and decomposition $W=U \oplus V$, we let
$r_{U,V}$ be the linear map which is $1$ on $U$ and $r$ on $V$.
Let us study how the symmetric space structure and the Jordan triple system 
of $M^\beta$ depends on $\beta$. In order to do this, we have to assume
that $p(o^+)=o^-$, i.e., that $\beta(o^+,o^-)=0$. In other words, 
$\beta$ is of the form
$$
\beta = \pmatrix{B_{1} & 0 \cr 0 & B_{2} \cr}, \quad
{\rm i.e.} \quad 
\beta((x,y),(x',y'))=\beta_1(x,x')+\beta_2(y,y') =
x^t B_1 x' + y^t B_2 y'
$$
with $x,x' \in \K^p$, $y,y' \in \K^q$ and
 matrices $B_1$ (size $p \times p$) and
$B_2$ (size $q \times q$) which are both either symmetric or 
skew-symmetric.
The property that $p(o^+)=o^-$ assures that the affine picture of
$p$ is a linear map from $M(q,p;\K)$ to $M(p,q;\K)$.
We calculate this map explicitly: let $X : \K^p \to \K^q$
be linear and consider its  graph
$$
\Gamma_X = \big\{ (x,X x)  \vert \, x \in \K^p \big\} \in \XX^+.
$$
Using that $B_1$ is invertible, we find that
the orthogonal complement of $\Gamma_X$ is given by
$$
\eqalign{
(\Gamma_X)^\beta & =
\{ (u,v) \in \K^p \times \K^q \vert \,
\forall x \in \K^p : \, \beta((u,v),(x,Xx))=0 \} \cr
&=
\{ (u,v) \in \K^p \times \K^q \vert \,
\forall x \in \K^p : \, u^t B_1x + v^t B_2 X  x   =0 \} \cr 
&=
\{ (u,v) \in \K^p \times \K^q \vert \,
\forall x \in \K^p : \, u^t x + v^t B_2 X B_1^{-1}  x   =0 \} \cr 
& =
\{ (u,v) \in \K^p \times \K^q \vert \,
\forall x \in \K^p : \, u = - B_1^{-1} X^t B_2 v  \} \cr
& =
\Gamma_{- B_1^{-1} X^t B_2} .
\cr}
$$
The linear picture of the map $p$ is therefore
$$
M(q,p;\K) \to M(p,q;\K), \quad
X \mapsto - B_1^{-1} X^t B_2 =
-(B_2 X B_1^{-1})^t = - ( \alpha X)^t
$$
with $\alpha(X)= B_2 X B_1^{-1}$, and the JTS of the symmetric space
$M^\beta$ is 
$$
T_\alpha(X,Y,Z) = X (\alpha Y)^t Z + Z (\alpha Y)^t X =  
X (B_2YB_1^{-1})^t Z + Z (B_2YB_1^{-1})^t X,
$$
which is the $\alpha$-homotope of $T(X,Y,Z)=X Y^t Z + Z Y^t X$. 

\ssk
Let us consider some special choices of $B_1$ and $B_2$.
Recall that $B_1$ has to be invertible, and $B_1$ and $B_2$ are either
both symmetric or both skew-symmetric.
 In case $\K=\R$ and $\K=\C$, topological connected components of
the symmetric spaces we obtain are homogeneous and of the
type $M^\beta \cong \OO(B_1 \oplus B_2)/(\OO(B_1) \times \OO(B_2))$.
The isomorphism classes of symmetric
spaces obtained for invertible matrices $B_2$ are listed in [Be00,
Tables XII.4]:

\ssk
\item{$\bullet$} Symmetric case: by choosing for $B_1$ and $B_2$ invertible
diagonal
matrices having coefficients  $b_{ii} \in \{ \pm 1 \}$, 
we get all symmetric spaces of the type
$\OO(l+j,k+i)/(\OO(l,k) \times \OO(i,j))$.
\item{$\bullet$} Skew-symmetric case: by choosing 
the usual matrices of symplectic forms (the dimensions have to be even),
we get all symmetric spaces of the type
$\Sp(l+j, \K)/(\Sp(l,\K) \times \Sp(j,\K))$.

\ssk \nin
Standard examples of deformations to degenerate cases are obtained
by taking some of the diagonal coefficients of $B_2$ equal to zero.
Since orthogonal groups of degenerate forms can be described as
semidirect products of vector groups, general linear groups and orthogonal
groups of non-degenerate forms, it is possible to give explicit descriptions
of such symmetric spaces. Note that in case $p$ even and $q$ odd a new
series of symmetric spaces arises which does  not arise as
 deformation of a semisimple series:
choose $B_1$ to be a symplectic form and $B_2$ skew-symmetric and non-zero.

\msk \nin {\bf 6.6. Other examples, and some final comments.}
Similar examples 
of deformations of semi-simple symmetric spaces
arise also in all other classical series of symmetric spaces (see the
classification in [Be00, Chapter XII]) and also for about half of 
the simple exceptional symmetric spaces (cf.\ comments in [Be00, Section 0.6]).
In order to keep this work in reasonable bounds, we cannot go here into
further details. Let us just mention that all classical groups arise
as group cases by our construction (the general linear groups are related
to associative algebras, see Section 0.1, and all orthogonal, symplectic
and unitary groups are related to certain Jordan triple systems); it
is then an interesting problem to determine the
deformations of these symmetric
spaces that are again of groupe type. In other words, we raise the problem:
{\it
which deformations of group type symmetric spaces come from deformations
of the overlying group structure}
(as was the case for the general linear group in Section 0.1)? 
This question should not be confused with another, also very natural
problem: 
{\it if $(M_\alpha)_ {\alpha \in \Svar(T)}$ is a deformation of the
symmetric space $M=M_{\id}$, in which sense are then the transvection
groups $G_\alpha:=G(M_\alpha)$ 
deformations of the transvection group $G=G(M)$?}
They will in general not be ``contractions of $G$" in the sense of
[DR85] since the dimension of $G(M_\alpha)$ will not be constant as 
a function of $\alpha$.
Replacing transvection groups by automorphism groups leads to a variant
of this problem. 
It seems that contractions of $\Aut(M)$ in the sense of [DR85] 
can only be obtained by using some additional parameter which in some sense
memorizes the structure of the stabilizer group of a point.

%Question: maybe in these group cases one can give an associative version of 
%the above associative construction  in terms
%of linear relations (which then would also allow to deform the group
%law, not only the symmetric spaces structure, in the group cases).

\sectionheadline
{Appendix A. Symmetric spaces}

\nin {\bf A.1. Reflection spaces with base point: a group-like concept.} 
A {\it reflection space} is a set $M$ together with a family
$(\sigma_x)_{x \in M}$ of {\it symmetries} $\sigma_x:M \to M$ 
such that

\ssk
\item{(S1)} $\sigma_x(x)=x$
\item{(S2)} $\sigma_x \circ \sigma_x = \id_M$
\item{(S3)} for all $x,y \in M$, $\sigma_x \sigma_y \sigma_x =
\sigma_{\sigma_x(y)}$.

\ssk \nin
The map
$$
\mu: M \times M \to M, \quad (x,y) \mapsto \mu(x,y):=\sigma_x(y)
$$
is called the {\it multiplication map}.
The concept of a reflection space is
``affine'', in the sense that it  does not refer to a distinguished base
point. On the other hand, for some purposes it is necessary to choose
a base point $o \in M$. 
Following Loos [Lo69], we then define the {\it quadratic map} and the
{\it square}, for $x \in M$,
$$
\QQ(x):=\QQ_o(x):= \sigma_x \circ \sigma_o: \, M \to M, \quad \quad
x^2:= \QQ(x)o=\sigma_x(o).
$$
The symmetry $\sigma_o$ is called {\it inversion} and is sometimes denoted by
$x^{-1}:=\sigma_o(x)$. Note that $\mu(x,y)$ can be recovered from
these data via 
$$
\mu(x,y)=\sigma_x (y)=\sigma_x \sigma_o \sigma_o (y) =
Q(x) j(y) = Q(x)Q(y)^{-1} y.
$$ 
The following relations are easily checked
(cf.\ [Lo69] or [Be00, Chapter I]): for all $x,y \in M$,

\ssk
\item{(SB1)} $\QQ(\QQ(x)y) = \QQ(x) \QQ(y) \QQ(x)$ 
\item{(SB2)} $\QQ(x^{-1}) = \QQ(x)^{-1}$
\item{(SB3)} $x^{-1} = \QQ(x)^{-1}x$

\Proposition A.2. There is an isomorphism of categories between the
category of reflection spaces with base point and the category of
 pointed sets $(M,o)$ equipped with two maps
$\QQ:M \times M \to M$, $(x,y) \mapsto \QQ(x)y$ and 
$\sigma_o:M \to M$, $x \mapsto
x^{-1}$ such that {\rm (SB1) -- (SB3)} hold. 

\Proof.
Given $(M,o;\QQ,\sigma_o)$, we define
$$
\mu(x,y):= \sigma_x (y):= \QQ(x)y^{-1} = \QQ(x)\QQ(y)^{-1}y,
$$
and (S1) -- (S3) are easily proved from (SB1) -- (SB3).
(For (S3), start by proving that $\sigma_o \QQ(x) \sigma_o = \QQ(x)^{-1}$.)
This construction clearly is inverse to the one described above, and
it is also clear that we get the same morphisms.
\qed

\nin 
The various spaces $(M,a,\QQ_a,\sigma_a)$ for $a \in M$ are in general
not isomorphic among each other 
(the map $\QQ(a)$ is an isomorphism from $(M,o)$ onto $(M,a^2)$ and not
onto $(M,a)$).
In order to change base points, 
we have the ``isotopy formula'':
$$
\QQ_a(x) = \sigma_x \sigma_a = \sigma_x \sigma_o (\sigma_a \sigma_o)^{-1}=
\QQ(x) \QQ(a)^{-1}.
$$

\msk \nin
{\bf A.3. Symmetric spaces.}
Assume $\K$ is a topological ring with dense unit group.
A {\it smooth symmetric space} is  a smooth manifold over $\K$ together with
a reflection space structure such that the multiplication map $\mu$
is smooth and

\ssk
\item{(S4)} 
the tangent map of the 
symmetry $\sigma_x$ at the fixed point $x$ is equal to the negative
of the identity of the tangent space $T_x M$.

\ssk 
\nin
(Cf.\ [Be06] for the notion of smpooth manifold over $\K$; if $\K=\R$ or
$\K=\C$, it coincides with the usual notion.)
We will use the term {\it symmetric space}, somewhat informally, also
in other, purely algebraic contexts, whenever the notion of tangent maps
of the symmetries $\sigma_x$ makes sense and (S4) is satisfied.
For instance, this will be the case if
there is some natural chart around $x$ such that $\sigma_x$ is,
in this chart, represented by the negative of the identity .

If $(M,o)$ is a smmoth symmetric space with base point, one may define its
{\it Lie triple system} in a similar way as the Lie algebra of a Lie group; 
this is done most naturally by using the quadratic maps $\QQ(x)$, cf.\
 [Lo69] for the real finite-dimensional and [Be06, Chapters 5 and  27]
for the general case.

Important examples of symmetric spaces arise from
the {\it group cases}: a group $G$ with the new product
$\mu(g,h)=\sigma_g(h)=gh^{-1}g$ satisfies (S1), (S2), (S3), and if $G$
is, e.g., a Lie group, then also (S4) holds.
Subspaces of $(G,\mu)$ 
are given by {\it spaces of symmetric elements}:
if $\tau$ is an involution of a group $G$ (anti-automorphism of order 2),
then $M:=\{g \in G | \tau(g)=g^{-1} \}$ is stable under $\mu$,
i.e., it is a subspace.

\sectionheadline
{Appendix B: Generalized projective geometries}

\nin The geometric counterpart of a Jordan pair is a {\it generalized
projective geometry}, and the one of a JTS a {\it generalized polar
geometry}. In this appendix we recall the basic definitions and refer
to [Be02] and [BeNe04] for further details.

\msk
\nin {\bf B.1. Affine pair geometries.}
A {\it pair geometry} is given by a pair of sets, denoted by
$(\XX^+,\XX^-)$, and a binary {\it transversality relation}
$(\XX^+ \times \XX^-)^\top \subset \XX^+ \times \XX^-$; 
we write $x \top \alpha$ or $\alpha \top x$ if $(x,\alpha)$ belongs to
$(\XX^+ \times \XX^-)^\top$, 
and we require that 
for all $\alpha \in \XX^\mp$ there exists $x \in \XX^\pm$ such
that $x \top \alpha$.
In other words, $\XX^\pm$ is covered by the sets
$$
\v_\alpha := \v_\alpha^\pm:= \alpha^\top :=
\{ x \in \XX^\pm  \vert \, x \top \alpha \}
\eqno{\rm (B.1)}
$$ 
as $\alpha$ runs through $\XX^\mp$.
An {\it affine pair geometry over a commutative ring $\K$} 
is a pair geometry such
that, for every  $\alpha \in \XX^\mp$,
 the set $\v_\alpha$ is equipped with the structure of an affine 
space over $\K$; then the sets $\v_\alpha$ will be called {\it affine
parts} or {\it affine charts} of $\XX^\pm$.
Given an affine pair geometry,
we let, for $(x,a), (y,a) \in (\XX^+ \times \XX^-)^\top$ and $r \in \K$,
$r_{x,a}(y):= r y$ be the product $r \cdot y$ in the $\K$-module
$\v_a$ with zero vector $x$, and we define the {\it structure maps} by
$$
\P_r:=\P_r^+: (\XX^+ \times \XX^- \times \XX^+)^\top  \to \XX^+, \quad
(x,a,y)   \mapsto \P_r(x,a,y) :=r_{x,a}(y),
\eqno{\rm (B.2)}
$$
where
$$
 (\XX^+ \times \XX^- \times \XX^+)^\top  := 
\{ (x,a,y) \in \XX^+ \times \XX^- \times \XX^+ \vert \,
x \top a, \, y \top a \},
\eqno{\rm (B.3)}
$$
and  dually the maps  $\P_r^-$ are defined.
Similarly, we get
structure maps defined by vector addition:
$$
\S:=\S^+: \,
 (\XX^+ \times \XX^- \times \XX^+ \times \XX^+)^\top  \to \XX^+, \quad
(x,a,y,z)   \mapsto \S(x,a,y,z) := y +_{x,a} z,
\eqno{\rm (B.4)}
$$
where the sum is taken in the $\K$-module $(\v_a,x)$ and
the set $(\XX^+ \times \XX^- \times \XX^+ \times \XX^+)^\top$ is defined
by a similar condition as in (B.3). Dually, $\S^-$ is defined.

\msk \nin {\bf B.2. Categorical notions.}
Assume $(\XX^+,\XX^-)$ and $(\YY^+,\YY^-)$ are affine pair geometries over
$\K$.

\ssk
\item{1.} Duality.
All axioms of an affine pair geometry 
appear together with their dual version: thus
$(\XX^-,\XX^+;\top;\P^-,\P^+,\S^-,\S^+)$ is again an affine pair geometry, called the
{\it dual geometry} of $(\XX^+,\XX^-;\top;\P^+,\P^-,\S^+,\S^-)$.
\item{2.} Homomorphisms. These are pairs of maps 
$(g^+:\XX^+ \to \YY^+,g^-:\XX^- \to \YY^-)$ which are compatible with
the transversality relations and with the structure maps, in the obvious 
sense.
\item{3.} Antiautomorphisms. These are isomorphisms $(g^+,g^-):(\XX^+,
\XX^-) \to (\XX^-,\XX^+)$ onto the dual geometry. If $g^-=(g^+)^{-1}$,
then $(g^+,g^-)$ is a  {\it correlation}, and if, moreover, there exists
a {\it non-isotropic point} (i.e., there is $x \in \XX^+$ with $x \top g^+(x)$),
then $(g^+,g^-)$ is called a {\it polarity}.
\item{4.} Adjoint (or: structural) pairs of morphisms.
These are pairs of maps $(f:\XX^+ \to \YY^+,h:\YY^- \to \XX^-)$ which
are compatible with transversality and with structure maps in the sense that
$x \top h(a)$ iff $f(x) \top a$ and
$$
f.\P_s(u,h.v,w)=\P_s(f.u,v,f.w), \quad
f.\S(u,h.v,w,h.z)=\S(f.u,v,h.w,z) .
$$
It is necessary to allow here also maps that are not everywhere defined 
(see Appendix C). A map $f:\XX^+ \to \XX^-$ is called {\it selfadjoint}
if $(f,f)$ is an adjoint pair of morphisms into the dual geometry
(compare with Section 1.12). 
\item{5.} Stability.
A pair geometry $(\XX^+,\XX^-,\top)$ is called {\it stable}
if any two points $x,y \in \XX^\pm$ are on a common affine chart, i.e., there
exists $a \in \XX^\mp$ such that $x,y \in \v_a$.

\msk \nin {\bf B.3. The affine pair geometry associated to a Jordan pair.}
Assume that $V=(V^+,V^-)$ is a linear Jordan pair over $\K$, which we
obtain in the
form $(V^+,V^-)=(\g_1,\g_{-1})$ for a 3-graded Lie algebra
$\g = \g_1 \oplus \g_0 \oplus \g_{-1}$ with $T^\pm(x,y,z)=[[x,y],z]$
(see Section 1.4).
Let $G:=\PE(V^+,V^-)=\langle U^+,U^+ \rangle$ be the 
{\it projective elementary group} associated to
this $3$-graded Lie algebra, i.e., the group generated by 
the vector groups $U^\pm = \exp(\g_{\pm 1})$,  let
$H$ be the intersection of $G$ with the automorphism group of the
$3$-grading (called the {\it (inner) structure group}) and let
$P^\pm := U^\pm H$. Since $U^\pm$ are vector groups, it is easily seen that 
$$
(\XX^+,\XX^-):=(\XX^+(V),\XX^-(V)):=(G/P^-,G/P^+)
$$
is an affine pair geometry. 
The base point $(eP^-,eP^+)$ is often denoted by $(o^+,o^-)$, and the
affine parts $(\v_{o^-},\v_{o^+})$ corresponding to this base point
are identified with $(V^+,V^-)=(\exp(\g^+).o^+,\exp(\g^-).o^-)$.
(In [BeNe04] a realization of this geometry and
of its transversality relation in
terms of {\it Lie algebra filtrations} is given.)

The construction is functorial:
a Jordan pair homomorphism induces a homomorphism of affine pair geometries
([Be02, Th.\ 10.1]).
In particular, an involution of a Jordan pair induces a polarity of
$(\XX^+,\XX^-)$. Since Jordan triple systems are nothing but Jordan pairs
with involution, we thus associate a polar geometry to any JTS.

The geometry associated to a Jordan pair $(V^+,V^-)$ is stable if and only if
the maps
$$
V^\mp \times V^\pm \to \XX^\pm, \quad (a,x) \mapsto \exp(a)\exp(x).o^\pm
$$
are surjective; then $\XX^\pm $ can be seen as 
$V^\pm \times V^\mp$ modulo the equivalence relation given by fibers of
this map (called {\it projective equivalence}; cf.\ [Lo94]).
Another equivalent version of stability is that the projective group $G$
admits the ``Harish-Chandra decomposition'' $G=\exp(V^+)\exp(V^-)H
\exp(V^+)$ and  $G=\exp(V^-)\exp(V^+)H\exp(V^-)$, where $H$ is the structure
group.

\msk \nin {\bf B.4.
Generalized projective geometries} are affine pair geometries
$(\XX^+,\XX^-)$ such that the structure maps satisfy some algebraic identities
(PG1) and (PG2) which should be seen as global analogs of certain
Jordan algebraic identies. 
Recall from B.1 that
$r_{x,a}(y)=\P_r(x,a,y)$, $r_{a,x}(b)=\P_r^-(a,x,b)$ (``left multiplication
operators"), and denote by $M^{(r)}_{x,y}(a):=\P_r(x,a,y)$,
$M^{(r)}_{a,b}(x):=\P_r(a,x,b)$ the ``middle multiplication operators".
Then the geometry $(\XX^+,\XX^-)$ associated to a Jordan pair $(V^+,V^-)$ 
satisfies the following properties ([Be02, Th.\ 10.1]):

\ssk
\item{(PG1)} For any invertible scalar $r \in \K$ and any transversal
pair $(x,a)$, the pair of
maps $(g,g')=(r_{x,a},r_{a,x}^{-1})$ is an automorphism of the affine pair
geometry $(\XX^+,\XX^-)$.
\item{(PG2)} For any $r \in \K$ and $a,b \in \XX^-$ with
$\v_a \cap \v_b \not= \emptyset$, the pair of middle multiplications
$(f,h)=(M^{(r)}_{a,b},M^{(r)}_{b,a})$ is an adjoint pair of (locally defined)
morphisms from $(\XX^+,\XX^-)$ to $(\XX^-,\XX^+)$, and dually.

\ssk \nin
Applied to elements of $\XX^+$, resp.\ $\XX^-$, (PG1) and (PG2) are indeed 
algebraic identities for
the structure maps. One should think of (PG1) as a global analog of
the defining identity (LJP2) (Section 1.4) of a Jordan pair, and of
(PG2) as the global analog of the fundamental formula.

\sectionheadline{Appendix C: Globalization of structural transformations}

\nin
{\bf C.1. Homomorphisms and structural transformations.}
There are two ways to turn Jordan pairs over $\K$ into a category:
one defines {\it homomorphisms} in the usual way as pairs of linear maps
$(h^+,h^-):(V^+,V^-) \to (W^+,W^-)$ which are compatible with the
maps $T^\pm$, and one defines a {\it structural transformation} to be a
pair of linear maps  $f:V^+ \to W^+$, $g:W^- \to V^-$ such that
$$
T^+_W(f(x),y,f(z)) = f (T^+_V(x,g(y),z)), \quad \quad
T^-_V (g(u),v,g(w)) = g (T^-_W(u,f(v),w))  
$$
(cf.\ [Lo94]).
In both cases we get a category; these categories are different, but
have essentially the same isomorphisms. 

\msk \nin {\bf C.2. Globalization: example of projective spaces.}
A homomorphism of Jordan pairs 
$(\phi^+,\phi^-):(V^+,V^-) \to (W^+,W^-)$ always lifts to a homomorphism
% $(\Phi^+,\Phi^-): (\XX^+(V),\XX^-(V)) \to(\XX^+(W),\XX^-(W))$ 
of associated geometries ([Be02, Th.\ 10.1]),
whereas for structural transformations the situation is more difficult.
In order to show that structural transformations do not always extend
globally,
consider the example of ordinary projective geometry $\XX^+ = \K \PPP^n$
 over a field. Its Jordan pair is $(W,W^*)$, with $W=\K^n$, and
$$
T^+(x,\lambda,y)=x \lambda(y) + y \lambda(x), \quad \quad
T^-(\lambda,x,\mu)=\lambda(x)\mu + \mu(x)\lambda .
$$
This Jordan pair is simple, and hence homomorphisms $(W,W^*)\to (V^+,V^-)$ 
are either injective or trivial. Clearly, every {\it injective} linear map
$\K^{n-1} \to \K^{m-1}$ of affine parts
lifts to a map $\K \PPP^n \to \K \PPP^m$ of associated projective spaces,
which shows that homomorphisms globalize, as stated above.
On the other hand, every linear map
$f:V \to W$ together with its dual map
defines a structural transformation $(f,g)=(f:V \to W,f^*:
W^* \to V^*)$. If $f$ is not injective, then $f$ does not lift
to a globally defined map of projective spaces.
The best we can do is to extend $f$ to a not everywhere 
defined map $\tilde f$ from $\K\PPP^n$ to $\K \PPP^m$.
The ``natural domain  of definition" of $\tilde f$ is
$$
D(\tilde f):= \Big\{ [ \pmatrix{v \cr r \cr} ] \, | \, v \in V, r \in \K :
\pmatrix{fv \cr r \cr}  \not= 0 \Big\} ,
$$
and then $\tilde f([ \pmatrix{v \cr r \cr} ])=[\pmatrix{fv \cr r \cr}]$.
The set $D(\tilde f)$
 is strictly bigger than the affine part $V$ (except for the trivial
case $f=0$).
 The purpose of this appendix is to define such a natural domain
of definition for any structural transformation of Jordan pairs.

\Lemma C.3.
If $(f,g)=(f:V^+ \to W^+,g:W^- \to V^-)$ is a structural transformation
between Jordan pairs, then
the pair $(V^+,W^-)$ with 
$$
\eqalign{
S^+:V^+ \times W^- \times V^+ \to V^-, & \quad (a,b,c) \mapsto
T^+_V (a,g(b),c) , \cr
S^-:W^- \times V^+ \times W^- \to W^-, & \quad (u,v,w) \mapsto
T^-_W (u,f(v),w)  \cr}
$$
is a linear Jordan pair, and the maps
$$
\eqalign{
\phi:=(\phi^+,\phi^-) &:=
(f,\id):(V^+,W^-) \to (W^+,W^-), \quad (x,y) \mapsto (f(x),y)  \cr
\psi:=(\psi^+,\psi^-) &:=
(\id,g):(V^+,W^-) \to (V^+,V^-), \quad (x,y) \mapsto (x,g(y))  \cr}
$$
are homomorphisms of Jordan pairs.

\Proof.
The first claim can be checked by a direct computation
which is completely analogous to the proof of Lemma 1.12.
%
%-- in fact, it may even be seen as a special case of this fact
%(consider the polarized JTS
%$$
%V = V^+ \oplus W^- \oplus V^- \oplus W^+
%$$
%(direct sum of the polarized JTS $V^+ \oplus V^-$ and $W^- \oplus W^+$).
%Then the map
%$$
%\alpha : V \to V, \quad (x,a,y,b) \mapsto (?)
%$$
%is structural and hence $V_\alpha$ is again a JTS.
%One checks that it is again polarized and hence corresponds to a Jordan
%pair. This Jordan pair is precisely $(V^+,W^-,S)$).
We show that $\phi$ is a homomorphism:
$$
\eqalign{
\phi^+ S^+(a,b,c)&=f T^+(a,gb,c) = T^+(fa,b,fc) = 
T^+(\phi^+ a,\phi^- b,\phi^+ c), \cr
\phi^- S^-(x,y,z)&= T^-(x,fy,z)= T^-(\phi^- x,\phi^+ y,\phi^- z). \cr}
$$
Similarly, we see that $\psi$ is a homomorphism.
\qed

\msk \nin
{\bf C.4. The globalization of a structural transformation.}
Assume $(f,g)$ is a structural transformation as above, and let
$\phi$ and $\psi$ the homomorphisms defined in the lemma.
By functoriality, we have induced homomorphisms on the level of geometries:
$$
\matrix{  &  (V^+,W^-) &   \cr
  {(f,\id) \atop }  \swarrow & &  \searrow {(\id,g) \atop } \cr
(W^+,W^-) &   & (V^+,V^-) \cr}
\quad \quad \quad \quad
\matrix{  &  \XX(V^+,W^-) &   \cr
  {(\Phi^+,\Phi^-) \atop }  \swarrow & &  \searrow {(\Psi^+,\Psi^-) \atop } \cr
\XX(W) &   & \XX(V) \cr}
$$
In general, if the first (resp.\ second) component of a Jordan pair 
homomorphism is injective, then so is the first (resp.\ second) component
of the induced map of 
geometries\footnote{$^1$}{\eightrm cf.\ [Lo94, 1.3 (3)] for the case of a stable
geometry; the proof of the corresponding fact for the general non-stable case 
is similar but rather lengthy and is therefore omitted here.}, and hence
 the first component $\Psi^+$ 
and the second component $\Phi^-$ are injective. 
The following maps together with their domains of definition are called
the {\it globalization of $(f,g)$}:
$$
\eqalign{
\tilde f: D(\tilde f) \to \XX^+(W), & \quad
 D(\tilde f):= \Psi^+(\XX^+(V^+,W^-)), \quad
\tilde f(\Psi^+(x)) := \Phi^+(x), \cr
\tilde g: D(\tilde g) \to \XX^-(V), & \quad
 D(\tilde g):= \Phi^-(\XX^-(V^+,W^-)), \quad
\tilde g(\Phi^-(a)) := \Psi^-(a). \cr}
$$
The following result will not be needed in this work and its proof is
therefore omitted:

\Theorem C.5. The pair $(\tilde f,\tilde g)$ is an adjoint
pair of morphisms between $\XX(V)$ and $\XX(W)$ whose domain of
definition $(D(\tilde f),D(\tilde g))$ is complete in the following
sense: whenever $x \in D(\tilde f)$, then $\tilde f(x)^\top \subset
D(\tilde g)$, and whenever $a \in D(\tilde g)$, then $\tilde g(a)^\top
\subset D(\tilde f)$.
Moreover,
$(D(\tilde f),D(\tilde g))$ is a generalized projective geometry, 
isomorphic to $\XX(V^+,W^-)$, with structure maps given by
$$
\P_r^+(x,a,y)=\P_r(x,\tilde g(a),y), \quad
\P_r^-(a,y,b)=\P_r^-(a,\tilde f(y),b),
\qeddis

\def\entries{

\[Be00 Bertram, W., {\it The Geometry of Jordan and Lie Structures},
Springer Lecture Notes in Mathematics vol. {\bf 1754},
Springer-Verlag, Berlin 2000

\[Be02 Bertram, W., ``Generalized projective geometries: General theory and equivalence with Jordan structures.''
    Advances in Geometry 2 (2002), 329-369

\[Be03 Bertram, W., ``The geometry of null systems, Jordan algebras and von Staudt's Theorem.''
     Ann. Inst. Fourier 53 (2003) fasc. 1, 193-225

\[Be06 Bertram, W., {\it Differential Geometry, Lie Groups and Symmetric Spaces 
over General Base Fields and Rings.}
 Memoirs of the AMS (to appear); 
      arXiv: math.DG/0502168 

\[BeHi00 Bertram, W. and J. Hilgert,
``Geometric Bergman and Hardy spaces'',
Michigan Math. J.\ {\bf  47} (2000), 235 -- 263 

\[BeHi01 Bertram, W. and J.\ Hilgert, ``Characterization of the Kantor-Koecher-Tits
algebra by a generalized Ahlfors operator'', J. of Lie Theory 11 (2001),
415--426

\[BeNe04 Bertram, W. and K.-H. Neeb, ``Projective completions of
Jordan pairs. I.'' J. Algebra {\bf 277}(2) (2004), 193 -- 225

\[BeNe05 Bertram, W. and K.-H. Neeb, ``Projective completions of
Jordan pairs. II'' Geom.\ Dedicata {\bf 112} (2005), 73 -- 113

\[DR85 Dooley, A.H. and J.W. Rice, ``On contractions of semisimple 
Lie groups'', Trans. A.M.S. {\bf 289} vol.1 (1985), 185 -- 202

%\[Fa Faraut, J., ``Int\'egrales de Riesz sur un espace 
%sym\'etrique ordonn\'e'', 

\[FaPe04 Faraut, J and M. Pevzner, ``Berezin kernels and Analysis on 
Makarevic spaces'', preprint; arXiv: math.RT/0411294

\[G64 Gerstenhaber, M., ``On the deformation of rings and algebras'',
Ann.\ Math.\ {\bf 79} (1964), 59 -- 103

%\[FK94 Faraut, J., and A. Kor\'anyi, {\it Analysis on Symmetric Cones}, 
%Clarendon Press, Oxford, 1994

%\[Hel80 Helgason, S., ``A Duality for Symmetric Spaces with Applications
%to Group Representations. III. Tangent Space Analysis''
%Adv. Math. {\bf 36} (1980), 297 -- 323

\[IW53 In\"onu, E., and E.\ P.\ Wigner, ``On the contraction of groups and
their representations'', Proc.\ Nat.\ Acad.\ Sci.\ USA {\bf 39} (1953),
510 -- 524

\[Lo69 Loos, O., {\it Symmetric Spaces I}, Benjamin, New York 1969

\[Lo75 Loos, O., {\it Jordan Pairs}, Springer LNM 460, Springer, New York 
1975

\[Lo77 Loos, O., {\it Bounded symmetric domains and Jordan pairs},
Lecture Notes, Irvine 1977

\[Lo94 Loos, O., ``Decomposition of projective spaces defined by unit-regular
Jordan pairs'', Comm.\ Alg.\ {\bf 22} (10) (1994), 3925 -- 3964

\[Ma73 Makarevic, B.O., ``Open symmetric orbits of reductive groups in
symmetric $R$-spaces'', Math.\ USSR Sbornik 20 (1973), no.\ 3, 406 -- 418

\[Ma79 Makarevic, B.O., ``Jordan algebras and orbits in symmetric 
$R$-spaces'', Trans.\ Moscow Math.\ Soc.\ {\bf 39} (1979)
(english translation, AMS, Trans.\ Moscow Math.\ Soc.\ 1981, issue 1,
169 -- 193)

\[McC04 McCrimmon, K., {\it A Taste of Jordan Algebras}, Springer 2004

\[Pe02 Pevzner, M., ``Analyse conforme sur les alg\`ebres de Jordan.''
[Conformal analysis on Jordan algebras]
J. Aust. Math. Soc. 73 (2002), no. 2, 279--299

\[Ri70 Rivillis, A.A., ``Homogeneous locally symmetric regions in homogeneous
spaces associated with semi-simple Jordan algebras'', Math.\ USSR Sbornik 11
(1970), no.\ 3 

%\[Sa Sabinin (reference to old literature on ``flocks'' ? cf. ..)
 
}

\references

\bye